\documentclass[twocolumn,5p]{elsarticle}
\usepackage{graphicx}%
\usepackage{multirow}%
\usepackage{amsmath,amssymb,amsfonts}%
\usepackage{amsthm}%
\usepackage{mathrsfs}%
\usepackage[title]{appendix}%
\usepackage{xcolor}%
\usepackage{textcomp}%
\usepackage{manyfoot}%
\usepackage{booktabs}%
\usepackage{algorithm}%
\usepackage{algorithmicx}%
\usepackage{algpseudocode}%
\usepackage{listings}%
\usepackage{subcaption}
\usepackage{float}
\usepackage{comment}

\usepackage{amssymb}
\usepackage{amsmath}
\journal{Biomedical Signal Processing and Control}

\begin{document}

\begin{frontmatter}



\title{The Effect of Prior Parameters on Standardized Kalman Filter-Based EEG Source Localization}


\author[TAU]{Dilshanie Prasikala}
\author[TAU]{Joonas Lahtinen}
\author[TAU]{Alexandra Koulouri}
\author[TAU]{Sampsa Pursiainen}
%% Author name

\affiliation[TAU]{organization={Mathematics Reseaerch Center, Computing Sciences, Tampere University},%Department and Organization
            addressline={Korkeakoulunkatu 3}, 
            city={Tampere},
            postcode={33720}, 
            state={Pirkanmaa},
            country={Finland}}

%% Abstract
\begin{abstract}

EEG Source localization is a critical tool in neuroscience, with applications ranging from epilepsy diagnosis to cognitive research. It involves solving an ill-posed inverse problem that lacks a unique solution unless constrained by prior knowledge. The Bayesian framework enables the incorporation of such knowledge, typically encoded through prior models. Various algorithms have been proposed for source localization, and they differ significantly in how prior knowledge is incorporated. Some approaches rely on anatomical or functional constraints, while others use statistical distributions or sampling-based techniques. In this landscape, the Standardized Kalman Filter (SKF) represents a dynamic Bayesian approach that integrates temporal modeling with a Gaussian prior structure. It addresses the depth bias, a common limitation in source localization, through a post-hoc standardization step that equalizes sensitivity across cortical depths and makes deep activity detection feasible.

This study focuses on the development and optimization of Gaussian prior models within the SKF framework for simultaneous cortical and sub-cortical activity detection. Synthetic data similar to the P20 / N20 component of the somatosensory evoked potentials (SEP) was used to identify effective prior parameter configurations for reconstructing both deep and superficial sources under different noise levels. We also investigated the role of RTS smoothing in enhancing source separability. Our results indicate that raising the standardization exponent to 1.25, along with smoothing, significantly improves depth localization accuracy at low noise levels.

By contributing to the relatively underexplored area of prior development in dynamic Bayesian frameworks, this work supports the viability of distributional estimation methods and provides a basis for future research into more flexible and adaptive prior structures for source localization.
\end{abstract}

\begin{graphicalabstract}
\end{graphicalabstract}

\begin{highlights}
\item This study introduces a systematic framework for tuning prior parameters of the Standardized Kalman Filter (SKF) to enhance spatio-temporal EEG source localization accuracy, especially for deep brain regions.
\item Our results suggest that SKF framework, when properly parameterized, mitigates depth bias and enables accurate tracking of spatio-temporal brain activity.
\end{highlights}

%% Keywords

\end{frontmatter}




\section{Introduction}\label{sec1}

Detecting deep brain activity is critical for understanding the functional organization of the brain and for diagnosing and treating neurological disorders such as Alzheimer’s disease \citep{laxton2013deep} and Parkinson’s disease \citep{obeso2008functional}. Non-invasive imaging modalities that provide both temporal and spatial resolution are essential for this task. Among these, electroencephalography (EEG) is a widely used technique due to its high temporal resolution, cost-effectiveness, and clinical accessibility. Recent studies have confirmed the ability of EEG to detect signals from subcortical regions, including the thalamus and hippocampus, particularly when advanced source reconstruction techniques are employed \citep{seeber2019subcortical}. However, accurately localizing deep sources remains a major challenge, primarily due to the attenuation of deep signals at the scalp and the ill-posed nature of the EEG inverse problem.

Many source localization algorithms solve each time step separately and ignore the temporal correlation of the data. Therefore, they can be referred to as instantaneous methods \cite{galka2004solution}. However, it can be observed that taking the temporal nature of data into account and including data from several points in time provides more information to the solution of the ill-posed inverse problem of source localization. Following this logic, the source localization problem can be viewed as a "dynamical" problem, and many algorithms have been developed to solve this dynamic problem \cite{somersalo2003non, schmitt2002efficient, brooks1999inverse, dannhauer2013spatio}.

Bayesian methodology provides a powerful framework for dynamic source estimation by modeling brain activity as a time-evolving process within a state-space system, where neural sources are latent variables updated recursively as new observations arrive. This approach reduces the computational demands of solving the full spatio-temporal inverse problem while allowing for real-time estimation and natural integration of prior models. To capture the non-linear and non-Gaussian nature of neural dynamics, nonlinear Bayesian filtering techniques such as particle filter-based source localization  \cite{sorrentino2009dynamical,sorrentino2010particle} has been developed. Recently, SESAME (SEquential Semi-Analytic Monte Carlo Estimator)  \citep{luria2024sesameeg,luria2020towards} has extended this line by combining sequential Bayesian updates with anatomical and sparsity-informed priors to robustly estimate both the number and location of active sources, particularly in deep brain regions.

While nonlinear Bayesian methodology offers this flexibility, it can be considered computationally demanding and often unsuitable for high-dimensional source spaces or real-time applications. To balance complexity and tractability, approximations such as the Extended Kalman Filter (EKF) \citep{ribeiro2004kalman} and Unscented Kalman Filter (UKF)  \citep{wan2000unscented} have been applied for EEG source localization. These variants rely on linearization or deterministic sampling and have been used with smoothness or anatomical priors in EEG source localization. Although these methods improve upon static solvers and can detect deeper activity to some extent, their effectiveness strongly depend on the accuracy of model assumptions and priors.

In scenarios where the assumption of linear Gaussian dynamics is reasonable, the standard Kalman filter provides an efficient and tractable approach. Within this framework, prior knowledge is typically encoded through simple dynamical models such as random walk or autoregressive priors, imposing temporal continuity in the estimated sources. Kalman filtering has been successfully applied to EEG source localization. For instance, \citep{galka2004solution} introduced a spatio-temporal Kalman filter using a first-order autoregressive prior with spatial coupling to capture structured source dynamics. \citep{barton2008evaluating} extended this by using a damped AR-2 prior to reflect oscillatory behavior, with parameters tuned via innovations likelihood. \citep{hamid2021source} applied a region-wise random-walk model with variable process-noise strength to enhance deep source detection. However, the depth bias introduced by the EEG lead field geometry remains a key limitation, often leading to underestimation of subcortical activity.

To mitigate this, the standardized Kalman Filter (SKF) algorithm was recently proposed \cite{Lahtinen2024SKF}. SKF introduces time-varying standardization weights, similar to those used in the sLORETA (standardized low-resolution brain
electromagnetic tomography) algorithm \citep{pascual2002standardized}, to adjust for the depth-dependent sensitivity of the EEG lead field. These weights are integrated into the recursive update step, effectively modifying the prior covariance structure to normalize the spatial influence of each source. While retaining the efficiency and recursive formulation of the standard Kalman filter, SKF significantly reduces depth bias and improves both spatial and temporal resolution in source reconstruction. Importantly, it maintains the assumption of Gaussian priors, but adapts their influence dynamically to favor deeper activity when supported by the data.

This paper builds on the SKF framework with the goal of advancing subcortical source localization using EEG. While SKF offers a promising solution, it currently lacks a principled strategy for the prior parameter selection, which critically affects the reconstruction accuracy as well as depth sensitivity. Furthermore, SKF is a computationally expensive algorithm, which is a common aspect in Kalman filtering. This makes careful prior parameter tuning essential to avoid excessive computation time and to ensure that the algorithm remains practical for use. So, as a first step, a systematic study of the prior parameter tuning will be carried out using synthetic data corresponding to somatosensory evoked potentials (SEP) occurring as a response to median nerve stimulation under three noise levels. In addition, we propose integrating a RTS smoother \citep{rauch1965maximum} on top of the tuned SKF algorithm and introducing a new tunable parameter called the "standardization exponent" as the exponent of the weighting matrix and investigating this parameter for values greater than 1, instead of the 0.5 used in the original SKF algorithm with the aim of enhancing the accurate recovery of deep neural activity with low dimensional source spaces.
 
\section{Theory}\label{sec2}

\subsection{Kalman Filtering Model}

Kalman filter is the closed-form solution to the Bayesian filtering equations, where the dynamic and measurement models are linear Gaussian. The state-space model is defined as follows: \cite{sarkka2023bayesian, simon2006optimal} 

\begin{equation}
\begin{aligned}
\mathbf{x}_t & =\mathbf{A}_{t-1} \mathbf{x}_{t-1}+\mathbf{q}_{t-1}, \\
\mathbf{y}_t & =\mathbf{L} \mathbf{x}_t+\mathbf{r}_t,
\end{aligned}
\end{equation}

where:
\begin{itemize}
    \item \( \mathbf{x}_t \in \mathbb{R}^n \) is the state vector (in this context, the current density),
    \item \( \mathbf{y}_t \in \mathbb{R}^m \) is the observation vector (electrical potentials measured outside the head),
    \item \( \mathbf{q}_{t-1} \sim \mathcal{N}(\mathbf{0}, \mathbf{Q}_{t-1}) \) is the process noise,
    \item \( \mathbf{r}_t \sim \mathcal{N}(\mathbf{0}, \mathbf{R}_t) \) is the measurement noise,
    \item \( \mathbf{A}_{t-1} \) is the state transition matrix,
    \item \( \mathbf{L} \) is the lead field matrix derived from the forward model,
    \item The prior distribution is \( \mathbf{x}_0 \sim \mathcal{N}(\mathbf{m}_0, \mathbf{P}_0) \), where \( \mathbf{m}_0 \) is the initial mean (assumed to be zero) and \( \mathbf{P}_0 \) is the initial covariance matrix.
\end{itemize}

The filtering model equations above are solved in closed form to obtain Gaussian distributions and Kalman filter prediction and update steps are then used to compute the parameters in the Gaussian distributions \citep{sarkka2023bayesian}. 

\subsection {Standardized Kalman Filtering Approach}

The SKF algorithm improves on the Kalman filtering algorithm through a standardization procedure equivalent to the one used in sLORETA \citep{pascual2002standardized}. However, it should be noted that Kalman Filtering is a time-dependent Bayesian estimation technique. Therefore, the standardization weights should be time-dependent unlike the sLORETA, which is a time-invariant source localization technique \citep{lahtinen2024standardized}. This is because, in the Kalman Filtering technique, the posterior distribution of the current state depends on the previous one.

Following the logic in standardization in sLORETA, given the states up to time step t-1, a dynamical standardization matrix $\mathbf{W}_t$ at each time step t for the MAP estimate of the Kalman Filter is introduced as follows(more on the derivation can be found in \cite{lahtinen2024standardized});

\begin{equation*}
\begin{aligned}
\mathbf{W}_t = \text{Diag} \left[ \mathbf{P}_{t|t-1}^{-1/2} \mathbf{K}_t \mathbf{S}_t \mathbf{K}_t^\top \mathbf{P}_{t|t-1}^{-1/2} \right]^{-\frac{1}{2}}
\end{aligned}
\tag{2}
\end{equation*}

where:
\begin{itemize}
  \item \( \mathbf{P}_{t|t-1} \) is the predictive covariance,
  \item \( \mathbf{K}_t \) is the Kalman gain,
  \item \( \mathbf{S}_t = \mathbf{L} \mathbf{P}_{t|t-1} \mathbf{L}^\top + \mathbf{R}_t \) is the innovation covariance.
\end{itemize}

Using the weighting matrix above, the standardized MAP estimate can be expressed as;

\begin{equation*}
\begin{aligned}
\mathbf{z}_{t \mid t} = \mathbf{W}_t \mathbf{P}_{t \mid t-1}^{-1 / 2} \mathbf{x}_{t \mid t}
\end{aligned}
\tag{3}
\end{equation*}

It has been proved \cite{lahtinen2024standardized} that the standardization weighting can be applied at the end of each time step, which gives more computation efficiency to the standardized Kalman filtering algorithm. Therefore, the above standardization can be used with the standard Kalman filtering algorithm to obtain the standardized Kalman filtering algorithm as follows; first, the prediction step of the algorithm is executed;

\begin{equation*}
\begin{aligned}
\mathbf{x}_{t \mid t-1} & = \mathbf{A}_t \mathbf{x}_{t-1 \mid t-1} \\
\mathbf{P}_{t \mid t-1} & = \mathbf{A}_t \mathbf{P}_{t-1 \mid t-1} \mathbf{A}_t^{\mathrm{T}} + \mathbf{Q}_t
\end{aligned}
\tag{4}
\end{equation*}

next, the update step is implemented;

\begin{equation*}
\begin{aligned}
\mathbf{S}_t & = \mathbf{L} \mathbf{P}_{t \mid t-1} \mathbf{L}^{\mathrm{T}} + \mathbf{R}_t \\
\mathbf{K}_t & = \mathbf{P}_{t \mid t-1} \mathbf{L}^{\mathrm{T}} \mathbf{S}_t^{-1} \\
\mathbf{x}_{t \mid t} & = \mathbf{x}_{t \mid t-1} + \mathbf{K}_t \left( \mathbf{y}_t - \mathbf{L} \mathbf{x}_{t \mid t-1} \right) \\
\mathbf{P}_{t \mid t} & = \mathbf{P}_{t \mid t-1} - \mathbf{K}_t \mathbf{S}_t \mathbf{K}_t^{\mathrm{T}} \\
\mathbf{W}_t & = \operatorname{Diag} \left( \mathbf{P}_{t \mid t-1}^{-1 / 2} \mathbf{K}_t \mathbf{S}_t \mathbf{K}_t^{\mathrm{T}} \mathbf{P}_{t \mid t-1}^{-1 / 2} \right)^{-1 / 2} \\
\mathbf{z}_{t \mid t} & = \mathbf{W}_t \mathbf{P}_{t \mid t-1}^{-1 / 2} \mathbf{x}_{t \mid t}
\end{aligned}
\tag{5}
\end{equation*}

In the current study, we suggest a tunable exponent \( \alpha \) on the weighting matrix, standardization exponent that can be applied to generalize the standardization:

\begin{equation*}
\begin{aligned}
\hat{\mathbf{z}}_{t|t}^{(\alpha)} = \mathbf{W}_t^{(\alpha)} \, \mathbf{P}_{t|t-1}^{-1/2} \hat{\mathbf{x}}_{t|t}
\end{aligned}
\tag{6}
\end{equation*}

where:

\begin{equation*}
\begin{aligned}
\mathbf{W}_t^{(\alpha)} = \text{Diag} \left[ \mathbf{P}_{t|t-1}^{-1/2} \mathbf{K}_t \mathbf{S}_t \mathbf{K}_t^\top \mathbf{P}_{t|t-1}^{-1/2} \right]^{-\alpha}
\end{aligned}
\tag{7}
\end{equation*}

\subsection {RTS Smoothing on the SKF algorithm}

In this study, we investigated the potential of using the Rauch, Tung, and Striebel \citep{rauch1965maximum} on the SKF algorithm to improve the quality of the reconstruction. This can be used to compute the closed-form smoothing solution to the Bayesian filtering equations. \citep{sarkka2023bayesian}. The main difference between the Kalman filter and the Kalman smoother is that the filter uses the whole set of measurements, while the filter uses only the measurements up to and including the current time step. Smoothers can be categorized into three groups and the RTS smoother discussed here falls into the category of fixed-interval smoothers where all the measurements in the interval are used to estimate the system state at any time in the interval \cite{grewal2014kalman}.

The process of smoothing is usually computationally demanding. However, the RTS smoother is considerably more efficient than other fixed-interval smoothers as it is not necessary to directly calculate the backward estimate or the covariance to arrive at the smoothed estimate \cite{grewal2014kalman}. 
Before the RTS smoother is initialized, the SKF algorithm is run from the beginning to the end of time.
Next, the RTS smoother is initialized as follows \cite{simon2006optimal};

\begin{equation*}
\begin{aligned}
\mathbf{x}_{T}^s & =\mathbf{x}_{T\mid T} \\
\mathbf{P}_{T}^s & =\mathbf{P}_{T\mid T}
\end{aligned}
\tag{3}
\end{equation*}

For $T = t-1, \ldots, 1, 0$, RTS Smoother equations can be executed as follows:

\begin{equation*}
\begin{aligned}
{\bf x}_{t+1\mid t} & = \mathbf{A}_t {\bf x}_{t\mid t} \\
\mathbf{P}_{t+1\mid t} & = \mathbf{A}_t \mathbf{P}_{t\mid t} \mathbf{A}_t^T + \mathbf{Q}_t \\
\mathbf{G}_t & = \mathbf{P}_{t\mid t} \mathbf{A}_t^T \left[\mathbf{P}_{t+1\mid t}\right]^{-1} \\
{\bf x}_{t}^s & = {\bf x}_{t\mid t} + \mathbf{G}_t \left[{\bf x}_{t+1}^s - {\bf x}_{t+1\mid t}\right] \\
\mathbf{P}_t^s & = \mathbf{P}_{t\mid t} + \mathbf{G}_t \left[\mathbf{P}_{t+1}^s - \mathbf{P}_{t+1\mid t}\right] \mathbf{G}_t^T
\end{aligned}
\tag{4}
\end{equation*}

In the current study, the effect of RTS smoothing was studied on the Standardized Kalman Filtering algorithm coupled with a suitable standardization exponent. Three values for the standardization exponent (1.00, 1.25, and 1.50) were explored to identify the suitable value to be applied.  

\section {Material and Methods}

\subsection {Model Parameters of the SKF Algorithm}

As a next step, the parameters needed to initialize and run the Standardized Kalman Filtering algorithm are considered. For the evolution of the SKF algorithm, a random walk model is assumed as there is no available physical model to describe the evolution of activity inside the brain. Here the current density is allowed to evolve in Gaussian distributed increments as it is assumed that $\mathbf{q}_{t}$ is normally distributed with zero mean and covariance $\mathbf{Q}_{t}$. Therefore, the transition matrix is taken as the identity matrix and this reduces the Standard Kalman transition model to,

\begin{equation*}
\begin{aligned}
\mathbf{x}_t &= \mathbf{x}_{t-1} + \mathbf{q}_{t},
\end{aligned}
\tag{8}
\end{equation*}

Here, the prior distribution is assumed to be Gaussian with covariance $\mathbf{P}_{0}$ and the prior covariance matrix is considered to be a diagonal matrix with equal entries $\boldsymbol{\theta}_0$. To select this value, the parameterization technique introduced in \cite{rezaei2020parametrizing} is used. Here, activity in individual source locations is assumed to be independent and identically distributed.

In the above-mentioned study, a new approach to the prior over-measurement signal-to-noise ratio (PM-SNR) is introduced. PM-SNR is referred to as the relative weight of the prior compared to the measurement noise. The PM-SNR parameter is a decibel value corresponding to the total prior standard deviation:

\begin{equation}
\mathrm{PM\text{-}SNR} = \mathrm{dB}(\sqrt{\boldsymbol{\theta}^{\text{\footnotesize (tot)}}_0}).
\end{equation}

It relates to the prior variance for each individual entry in the source space as \citep{rezaei2020parametrizing}:

\begin{equation} 
\boldsymbol{\theta}_0 = \frac{\boldsymbol{\theta}_0^{\text{\footnotesize (tot)}} \, \boldsymbol{\sigma}^2 \, \mathbf{A}^2}{\mathbf{N}},
\end{equation} 

where $\mathbf{A}$ is the measurement amplitude, $\boldsymbol{\sigma}$ is the noise standard deviation, and $\mathbf{N}$ is the number of sources in the source space. At the reference level $\mathrm{PM\text{-}SNR} = 0$ dB, the total standard deviation approximates the maximum amplitude of the spatiotemporal brain activity: $\sqrt{\boldsymbol{\theta}^{\text{\footnotesize (tot)}}_0} \approx \max |\vec{J}^p|$

Following a similar idea, evolution (process noise) can also be considered as an independent random variable $\mathbf{q}$ and the evolution prior (process noise) covariance matrix ($\mathbf{Q}_{t}$) a diagonal matrix with equal variances. Following from the additivity of variance under summation of Gaussian random variables, and because the prior is a Gaussian one, we can define a net evolution prior variance $\boldsymbol{\tau}^2$ as follows:

\begin{equation*}
\begin{aligned}
\boldsymbol{\tau}^2 = \sum_{i = 1}^K \boldsymbol{\tau}_i^2.
\end{aligned}
\tag{9}
\end{equation*}

Here, $\boldsymbol{\tau}_i$ is the evolution variance per time step and $K$ is the number of the time steps in the Markov chain. To keep the prior strength unchanged if the number of time steps is changed, we relate $\boldsymbol{\tau}^2$ to the prior variance as follows:

\begin{equation*}
\begin{aligned}
\boldsymbol{\tau}^2 &= \boldsymbol{\kappa}^2 \, \boldsymbol{\theta}_0 \\
\boldsymbol{\tau}_i^2 &= \boldsymbol{\kappa}^2 \, \frac{\boldsymbol{\theta}_0}{K}
\end{aligned}
\tag{10}
\end{equation*}

with $\boldsymbol{\kappa}$ corresponding to evolution prior signal-to-noise ratio: $\mathrm{EP\text{-}SNR} = \mathrm{dB}(\boldsymbol{\kappa})$.

Whenever the evolution drive is greater than zero, the system tries to advance the temporal progression more than what is suggested by the {\em a priori} knowledge of the signal amplitude. When it is less than zero, the evolution is weaker than expected {\em a priori}. That is, if the evolution drive is kept at a low value, the algorithm tends to be biased toward the source it first localizes.

\subsection{Simulation Setup}

\begin{figure*}[h!]
    \centering
    \begin{minipage}{0.45\textwidth}
        \centering
        \includegraphics[width=0.55\textwidth]{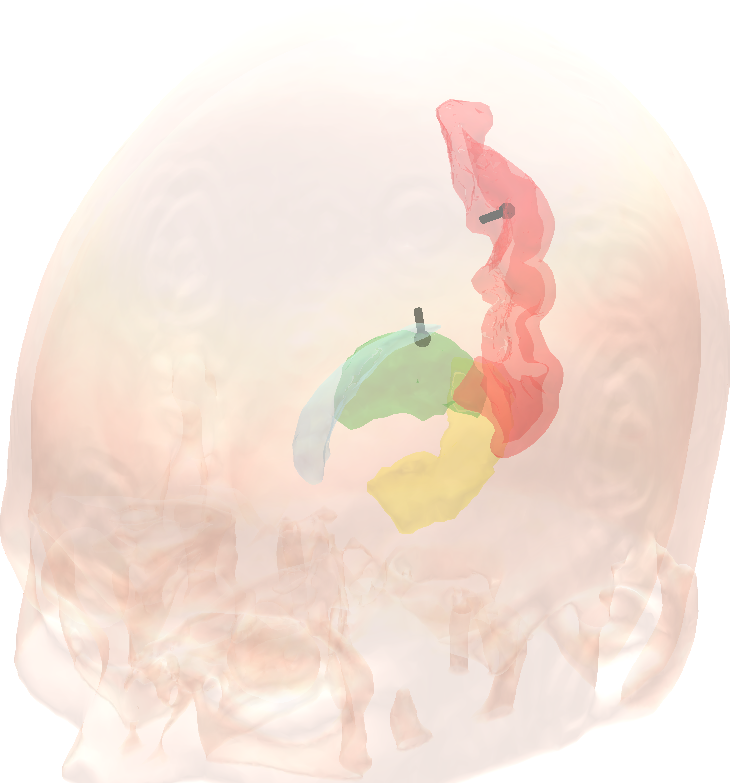} \\ \vskip0.1cm  \includegraphics[width=0.9\linewidth,trim=0 50 0 250, clip]{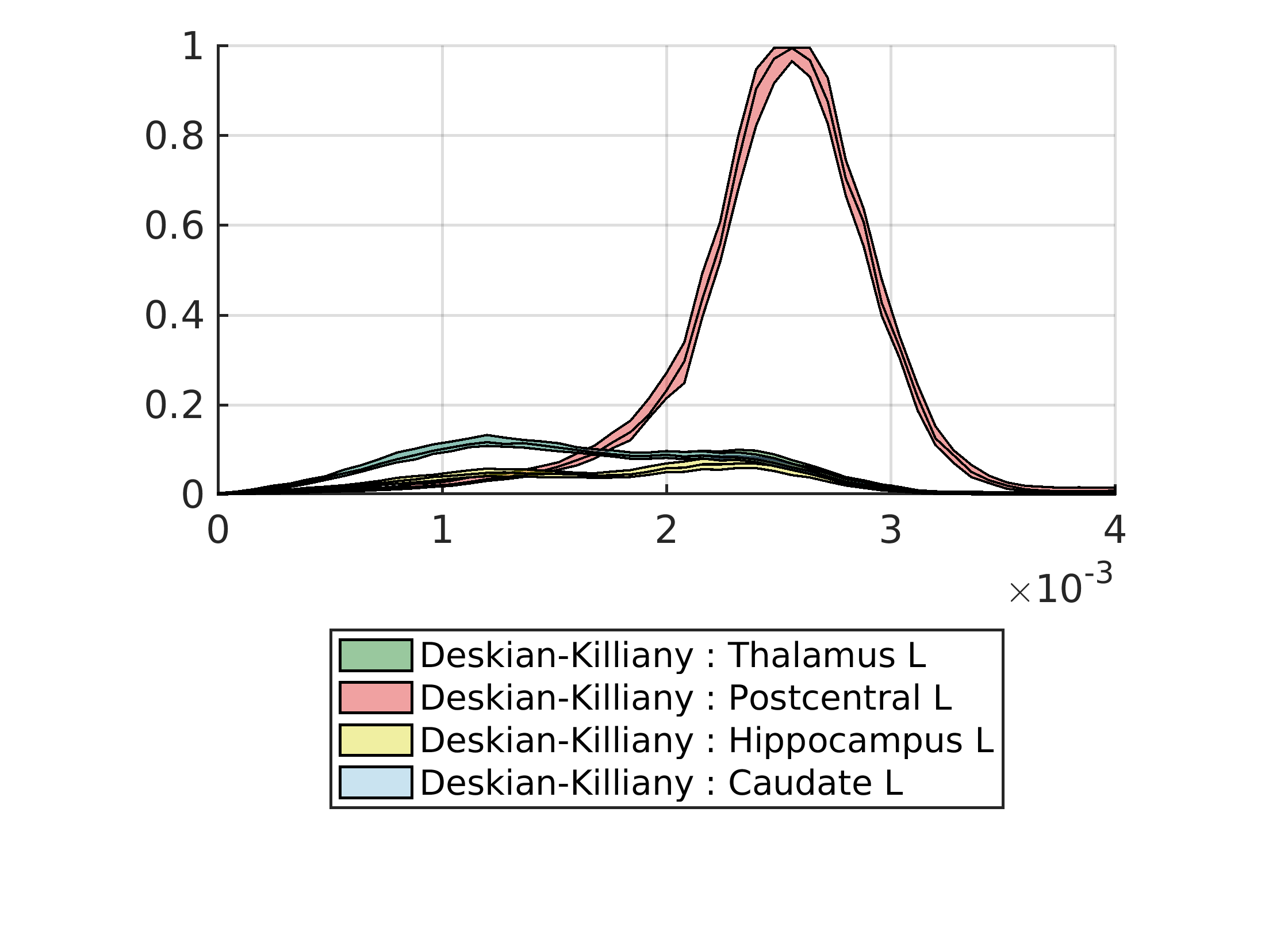} \\ (a)
    \end{minipage}
    \hfill
    \begin{minipage}{0.45\textwidth}
        \centering
        \includegraphics[width=0.8\textwidth]{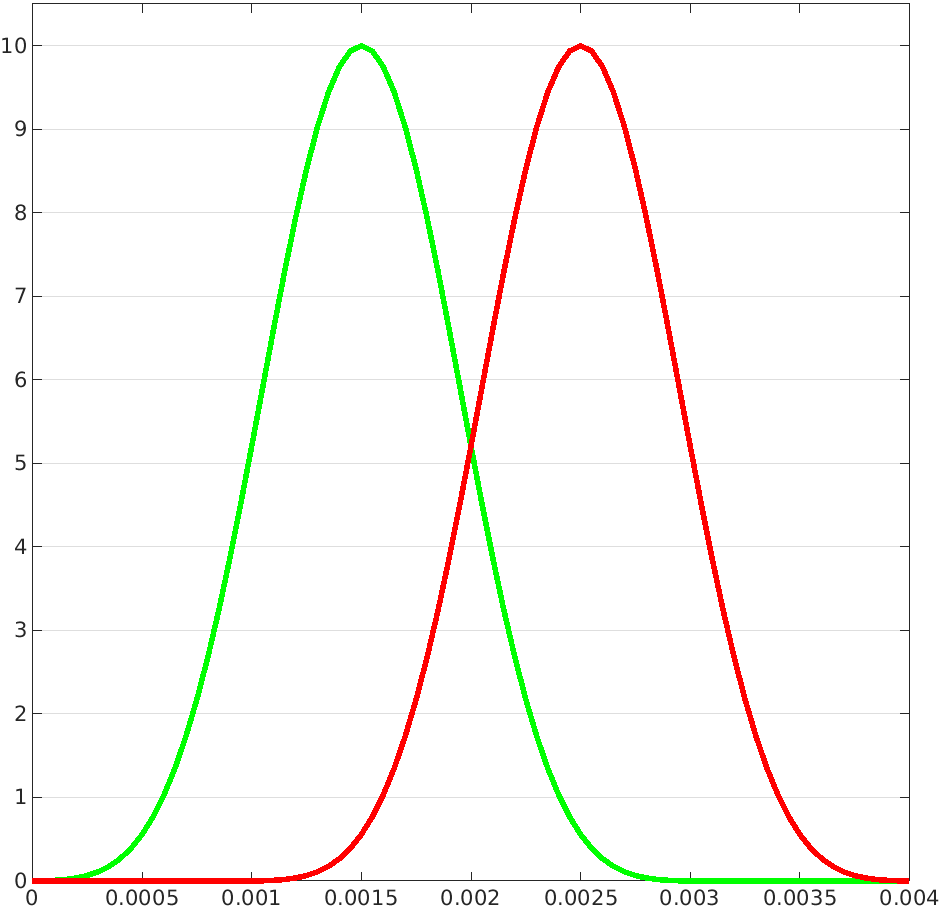} \\
        (b)
    \end{minipage}
    
    \caption{(a) - Source placement corresponding to P20/N20 (20 ms post-stimulus response) component of median nerve SEP, where the cortical component
is visible at the posterior wall of the central sulcus corresponding to a tangentially 
oriented source and the sub-cortical component is visible at the ventral posterolateral 
part of the thalamus. (b) - Time series plot for 4 ms simulated signal corresponding to 2 ms  activation of subcortical activity (Green curve) followed by 2 ms peak of cortical activity (Red curve) with a 2 ms overlap between the peaks. The x-axis represents the time in seconds and the y-axis represents the amplitude in nanoampere-meters.}
    \label{fig:source_placement}
\end{figure*}

\subsubsection{Head Model}

A multi-compartment head model was created using the ICBM152 2023b \citep{mazziotta2001probabilistic} template, which is openly available in Brainstorm \citep{tadel2011brainstorm}. The Desikan-Killiany atlas \cite{desikan2006automated} which is available in Brainstorm (BST) was used to visualize the cortical structures of the brain. Zeffiro Interface (ZI) \cite{he2020zeffiro} was used as the primary tool for analysis in this study. The head model created in BST was imported to the ZI  using the BST-2-ZI pipeline available in ZI's repository. 

\subsubsection{Lead Field Matrix}

To solve the forward problem, a boundary-fitted FEM (Finite Element Method) mesh with 3 mm resolution \cite{galaz2023multi} was generated using ZI. An EEG lead field matrix of 3000 source positions uniformly distributed in the active compartments was created through the divergence confirming source model \citep{pursiainen2016electroencephalography}. This lead field size was found to be sufficient to represent all the active compartments, and a larger size would have used more computational capacity.

\subsubsection{Simulated Sources}

In this study, two synthetic sources were placed in the head model; one in the deep part of the brain; specifically in the left ventral posterolateral thalamus, and the other in the posterior wall of the left central sulcus positioned to match the deep and superficial components of the P20/N20 (20 ms post-stimulus response) component of the somatosensory evoked potentials (SEP) occuring as a response to median nerve stimulation  \citep{buchner1994source,rezaei2021reconstructing}.The orientations of the sources were also adjusted to match the SEP components, where the superficial source was kept tangential (Figure \ref{fig:source_placement} (a)). In order to approximately match the SEP response, where the afferent volley propagates from the thalamus to the cerebral cortex, the signal was generated as a combination of two smooth overlapping  2 ms peak signals with a total duration of 4 ms, a 2 ms total overlap, and a 2 ms separation between peaks (Figure \ref{fig:source_placement} (b)).

\section{Results}

\subsection{Parameter Tuning at Different Noise Levels}

In the first stage of the study, the tuning of prior parameters was considered adding a varying amount of Gaussian zero-mean white noise into the synthetic measurements. Reconstructions were evaluated for the following three measurement noise peak-SNR levels: 10 dB, 20 dB, and 30 dB

At the highest noise level (10 dB), the best reconstruction quality was observed when both the EP-SNR and the PM-SNR were set to zero. Increasing EP-SNR or PM-SNR beyond zero degraded the quality of the reconstruction, leading to loss of spatial accuracy and shape of the signal. Furthermore, for the same magnitude of increment, changes in EP-SNR exhibited a more substantial negative effect compared to changes in PM-SNR (Figure \ref{fig:comparison_noise10}).

When the noise level was decreased to 20 dB, a similar pattern was observed. The best reconstruction continued to occur at low values of EP-SNR and PM-SNR (both at zero). However, the sensitivity of the reconstruction quality to changes in EP-SNR and PM-SNR was significantly lower compared to the 10 dB case, indicating greater robustness at moderate noise levels (Figure \ref{fig:comparison_noise20}).

At the lowest noise level (30 dB), the behavior of the system changed. Here, the best reconstructions were obtained at EP-SNR = 20 and PM-SNR = 0, or vice versa. Further increases beyond these values again lead to degradation of the source localization performance (Figure \ref{fig:comparison_noise30}).

Across all noise levels, even the best reconstructions exhibited a characteristic artifact: the presence of an echo of the deep activity around the time of the superficial peak. This echo-effect persisted regardless of parameter tuning, indicating a structural challenge in separating deep and superficial sources under noise.

\subsection{Effect of Smoothing on Source Reconstruction}

The effect of Rauch-Tung-Striebel (RTS) smoothing on the quality of the source reconstruction was evaluated. In general, smoothing was observed to suppress deep source activity under all conditions. However, adjusting the standardization exponent within the SKF algorithm modulated the effect of smoothing. Specifically, increasing the standardization exponent above 1.0 improved the reconstructions by reducing the influence of spurious deep activity and reducing the secondary peak at the time of superficial activation.

\begin{figure*}[ht!]
    \centering
    \begin{subfigure}[t]{0.32\textwidth}
        \centering    
        \includegraphics[width=1.08\linewidth,trim=0 125 0 0, clip]{figures/reconstructions/sta_exp_1.png}
        \caption{}
        \label{fig:fig1}
    \end{subfigure}
    \hfill
    \begin{subfigure}[t]{0.32\textwidth}
        \centering
        \includegraphics[width=\linewidth, trim=0 125 0 0, clip]{figures/results_1.25_smoothed/results_plot_inversion_data_kalman_prior_30_20_0.png}
        \caption{}
        \label{fig:fig2}
    \end{subfigure}
    \hfill
    \begin{subfigure}[t]{0.32\textwidth}
        \centering
        \includegraphics[width=\linewidth,trim=0 125 0 0, clip]{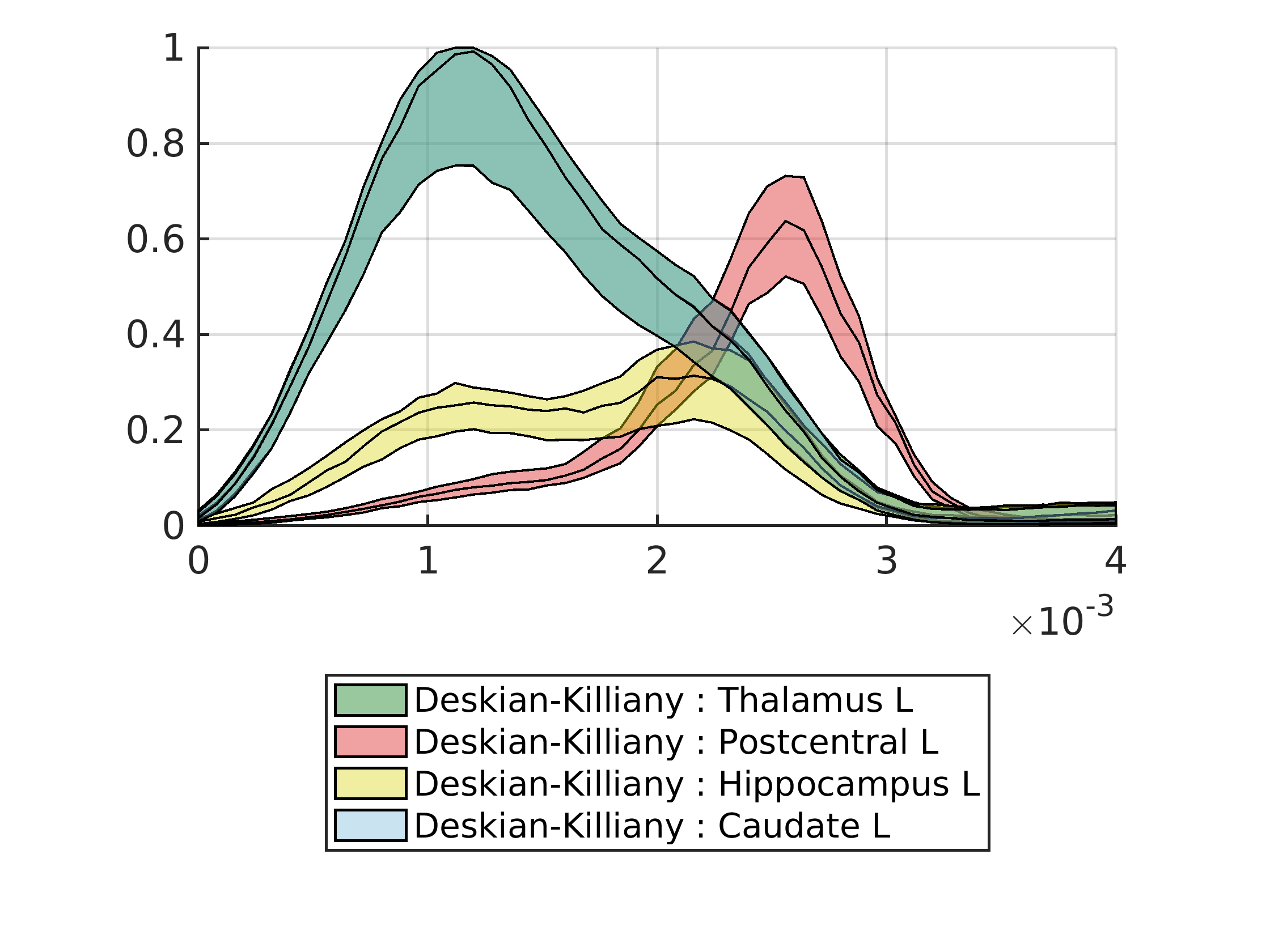}
        \caption{}
        \label{fig:fig3}
    \end{subfigure}
       
    \caption{Time series plots illustrating the effect of standardization exponent on RTS smoother for the SKF algorithm. standardization exponents 1.00 (left), 1.25 (middle) and 1.5 (right).The smoothing results are for the 30 db noise case with tuned EP-SNR and PM-SNR parameters(20 and 0 respectively). The colors of the curves are as in Figure \ref{fig:source_placement}.}
    \label{fig:standardization_exponent}
\end{figure*}

Three values for the standardization exponent (1, 1.25 and 1.5) were considered. It was found that a standardization exponent value of 1.25 yielded the overall best reconstructions. The value of 1 suppressed the deep peak, while the value of to 1.5 tended to overemphasize the deep peak, degrading the localization performance (Figure \ref{fig:standardization_exponent}).

\begin{figure*}[h!]
\centering
\noindent
\begin{minipage}[t]{0.46\textwidth}
    \centering    

    \begin{subfigure}{0.48\textwidth}
        \includegraphics[width=\linewidth, trim=0 135 0 0, clip]{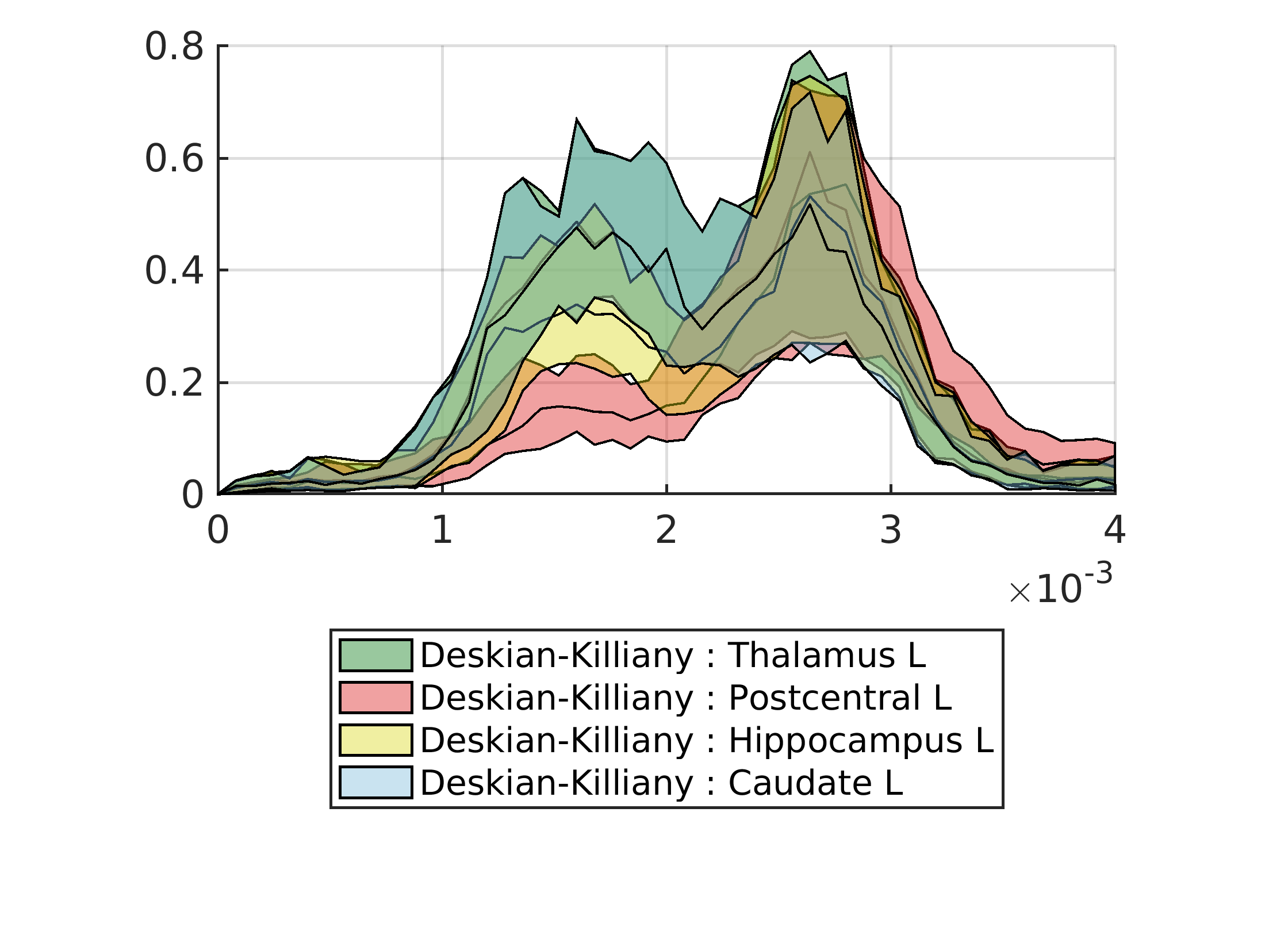}
        \caption*{EP 0, PM 0}
    \end{subfigure}
    \hfill
    \begin{subfigure}{0.48\textwidth}
        \includegraphics[width=\linewidth, trim=0 135 0 0, clip]{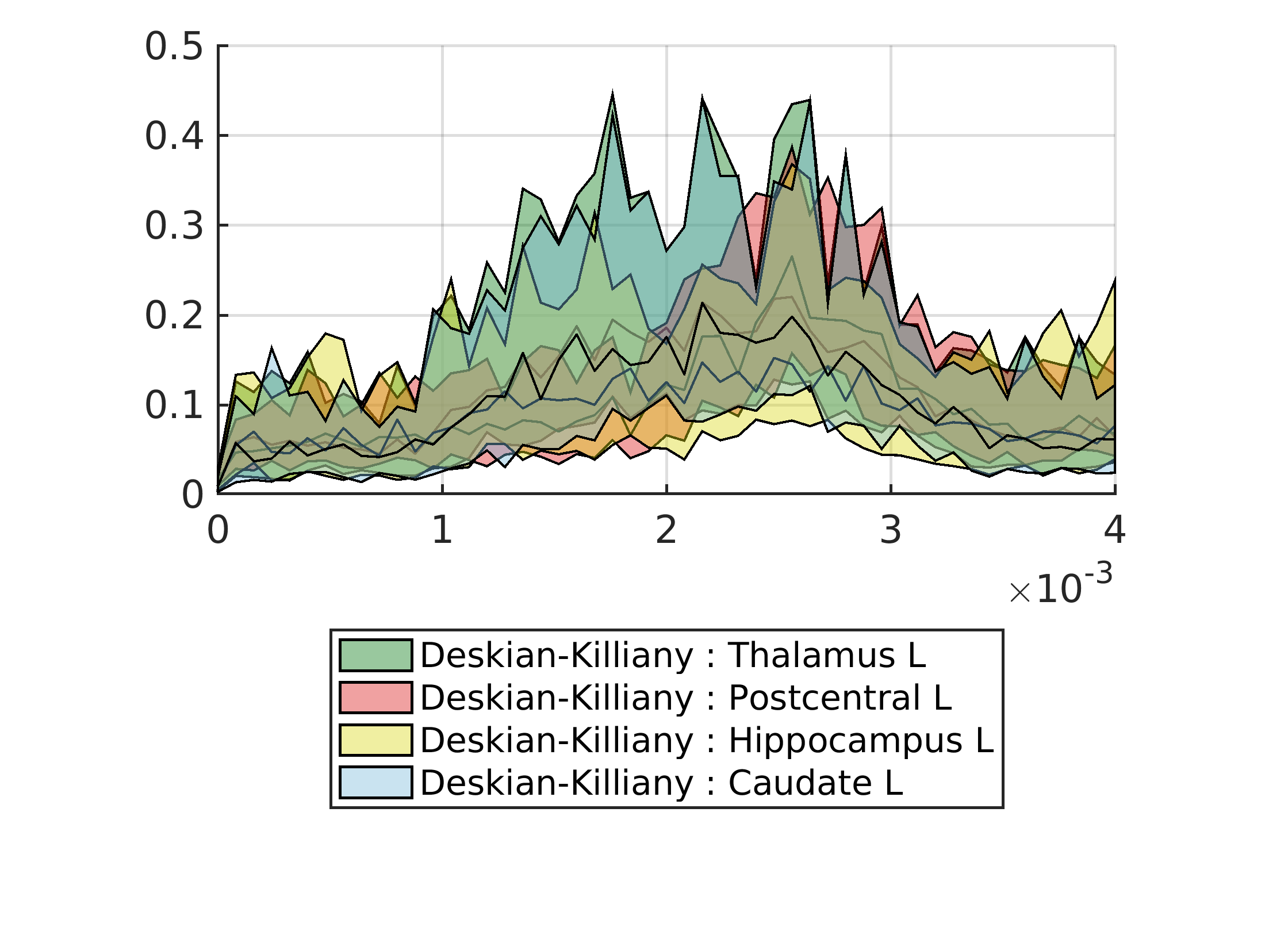}
        \caption*{EP 20, PM 0}
    \end{subfigure}

    \vspace{0.5em}

    \begin{subfigure}{0.48\textwidth}
        \includegraphics[width=\linewidth, trim=0 135 0 0, clip]{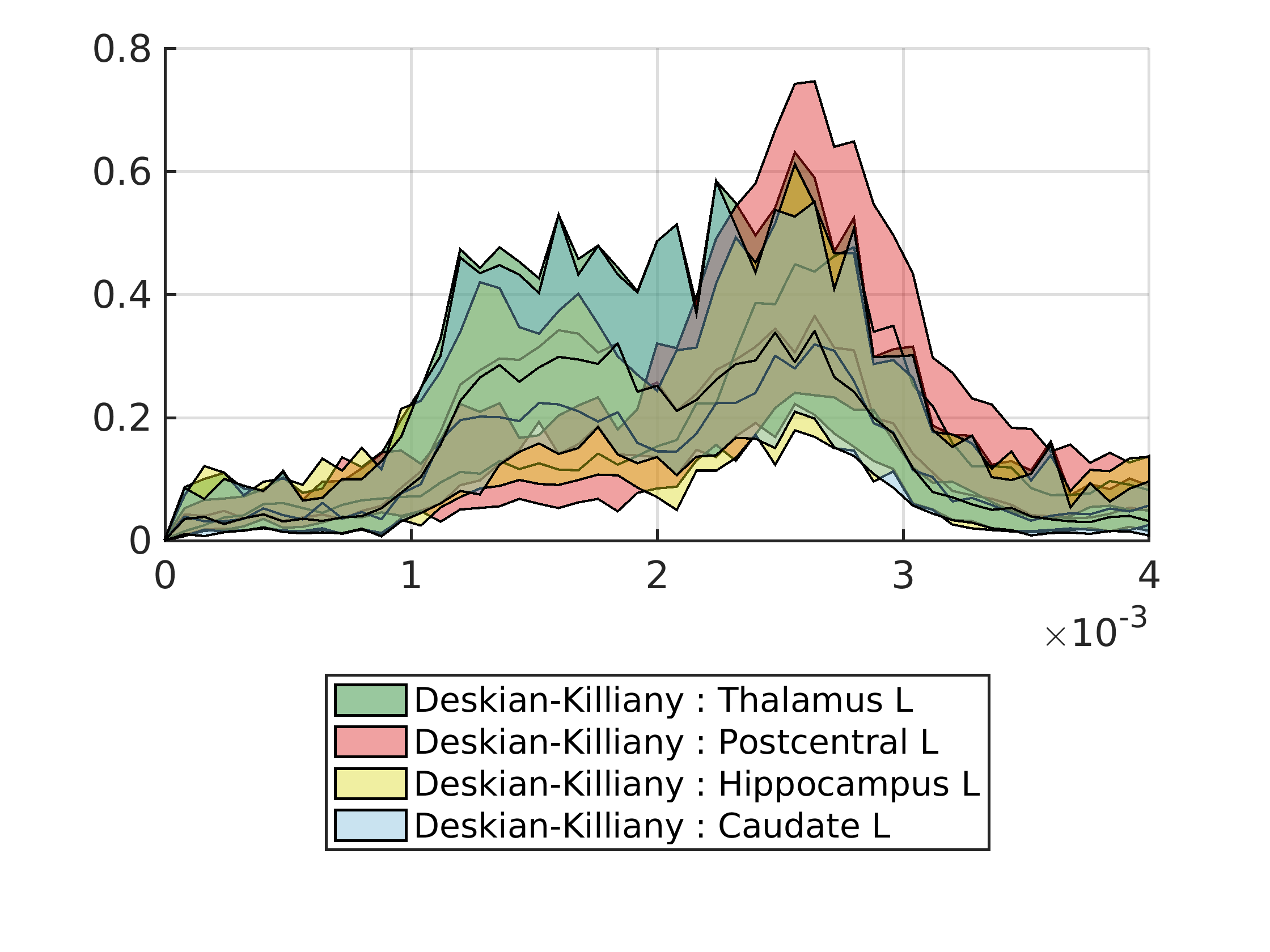}
        \caption*{EP 0, PM 20}
    \end{subfigure}
    \hfill
    \begin{subfigure}{0.48\textwidth}
        \includegraphics[width=\linewidth, trim=0 135 0 0, clip]{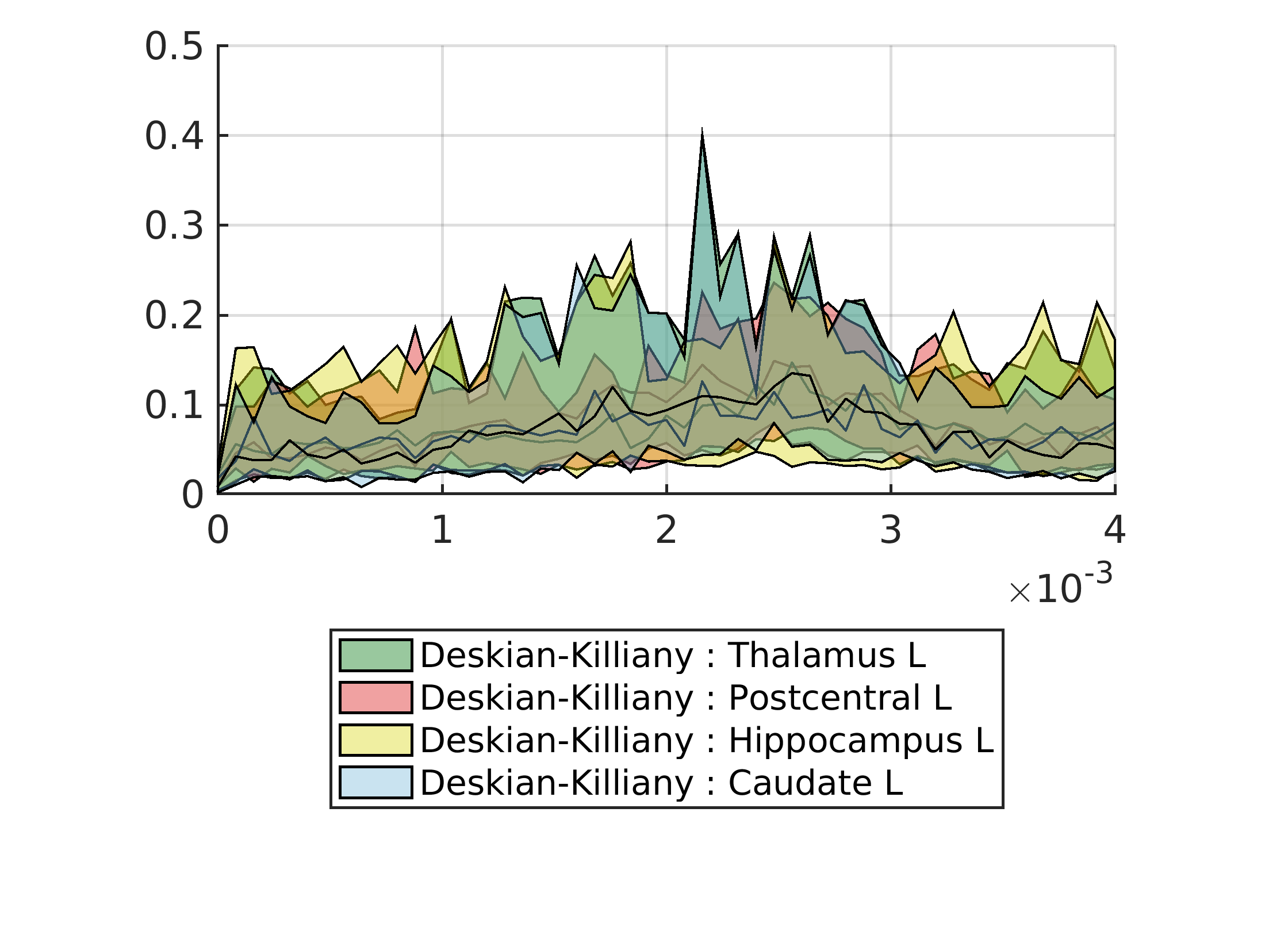}
        \caption*{EP 20, PM 20}
    \end{subfigure}
\end{minipage}
\hfill
\begin{minipage}[t]{0.46\textwidth}
    \centering

    \begin{subfigure}{0.48\textwidth}
        \includegraphics[width=\linewidth, trim=0 135 0 0, clip]{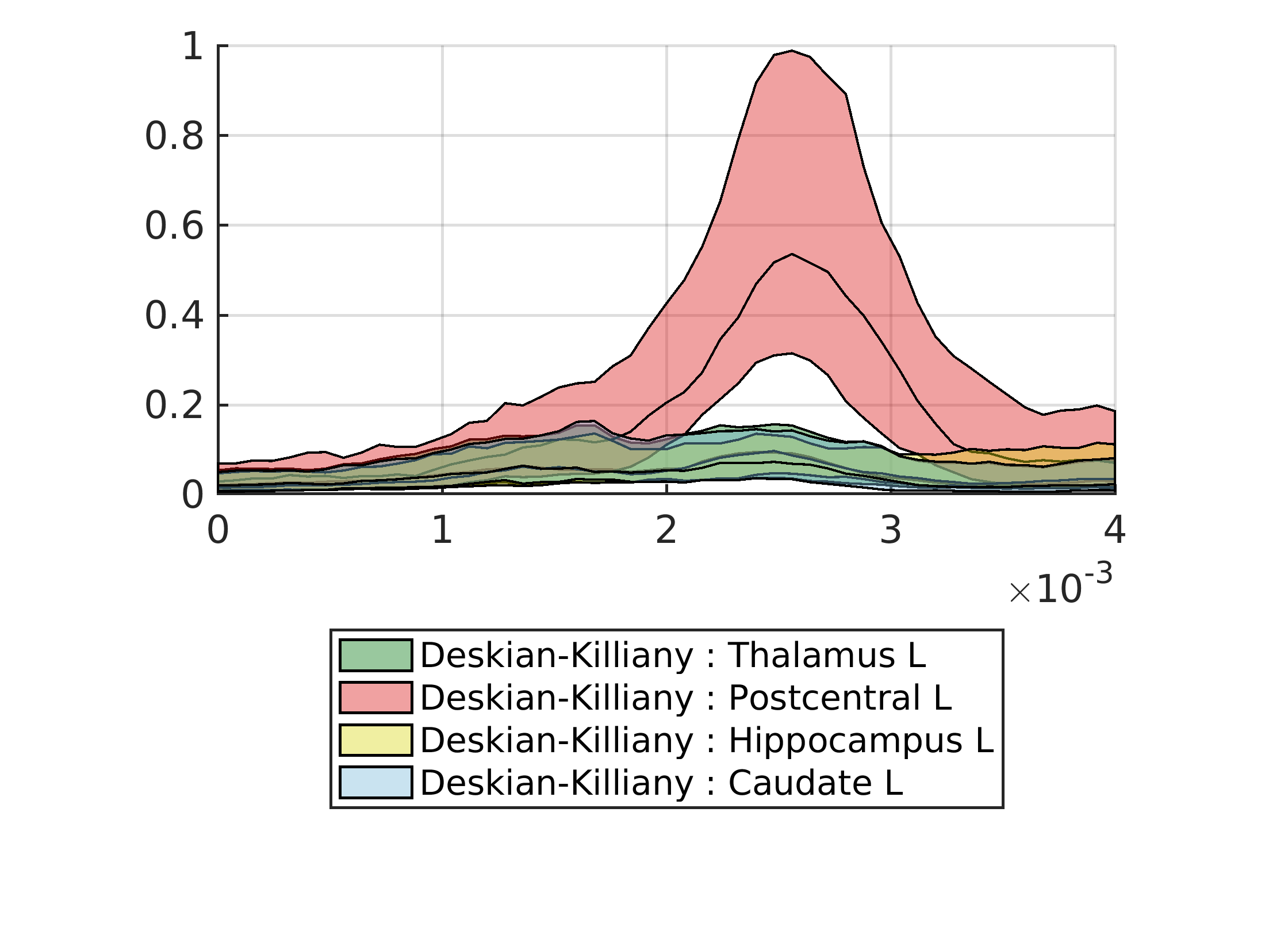}
        \caption*{EP 0, PM 0}
    \end{subfigure}
    \hfill
    \begin{subfigure}{0.48\textwidth}
        \includegraphics[width=\linewidth, trim=0 135 0 0, clip]{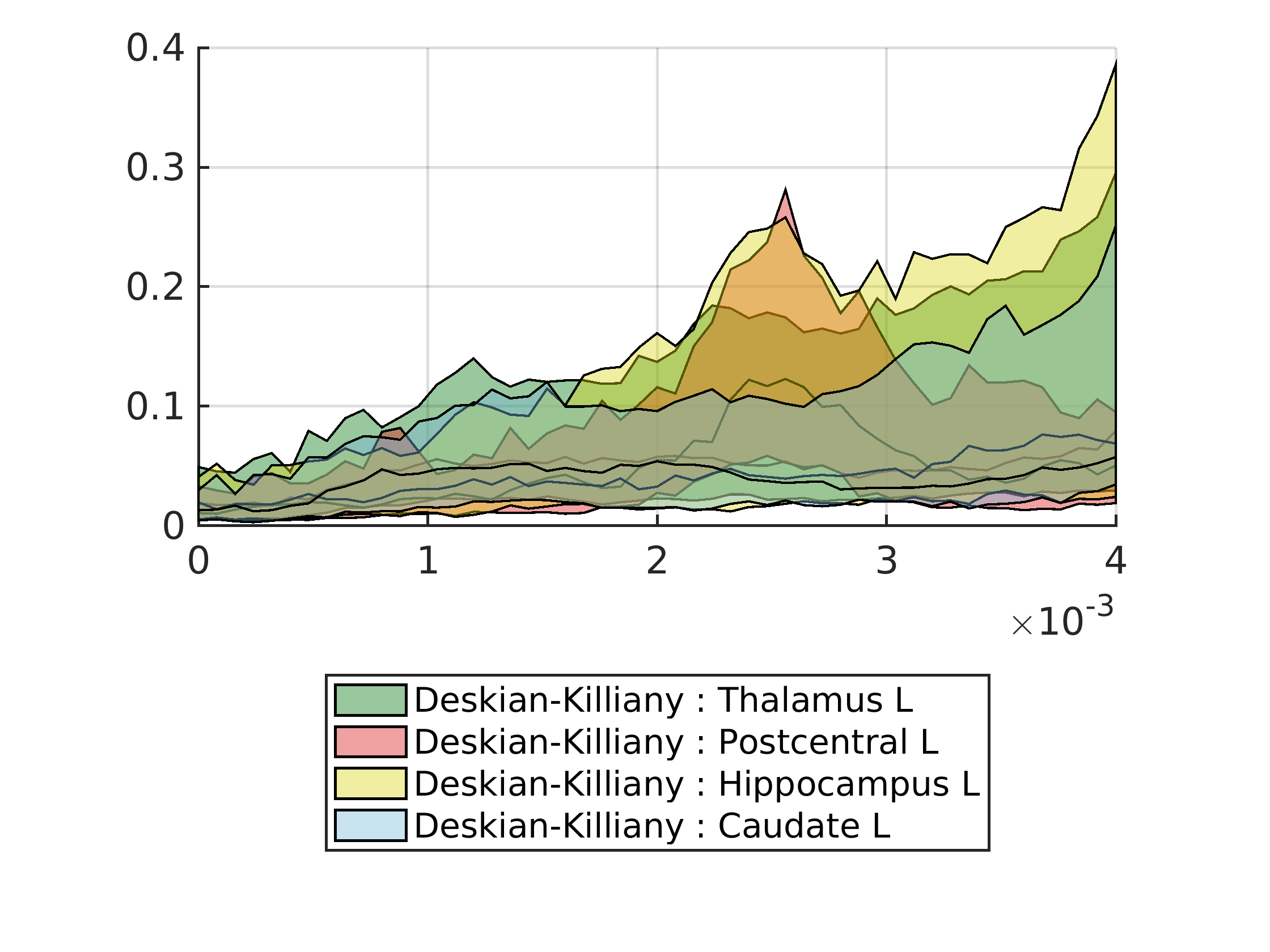}
        \caption*{EP 20, PM 0}
    \end{subfigure}

    \vspace{0.5em}

    \begin{subfigure}{0.48\textwidth}
        \includegraphics[width=\linewidth, trim=0 135 0 0, clip]{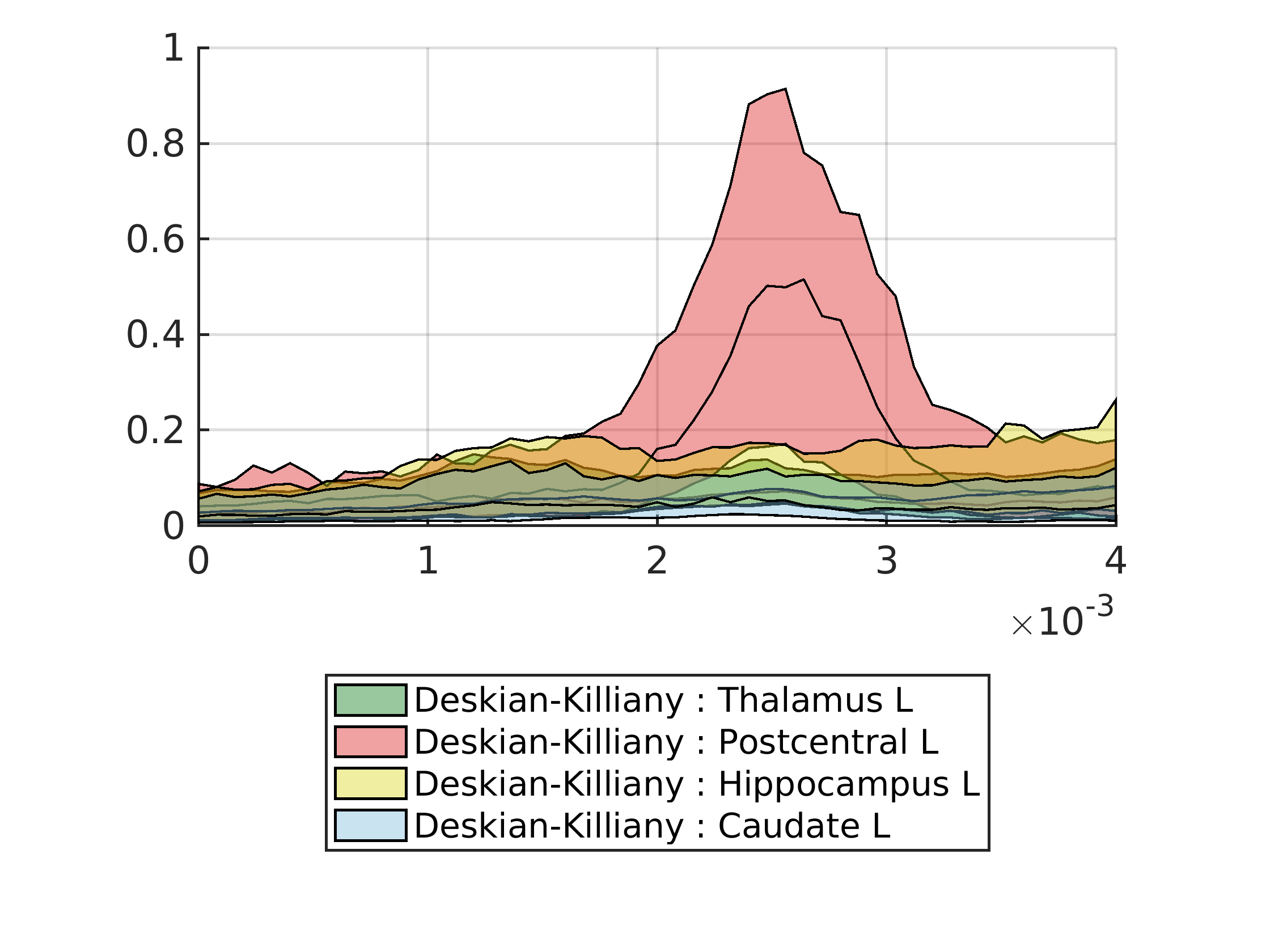}
        \caption*{EP 0, PM 20}
    \end{subfigure}
    \hfill
    \begin{subfigure}{0.48\textwidth}
        \includegraphics[width=\linewidth, trim=0 135 0 0, clip]{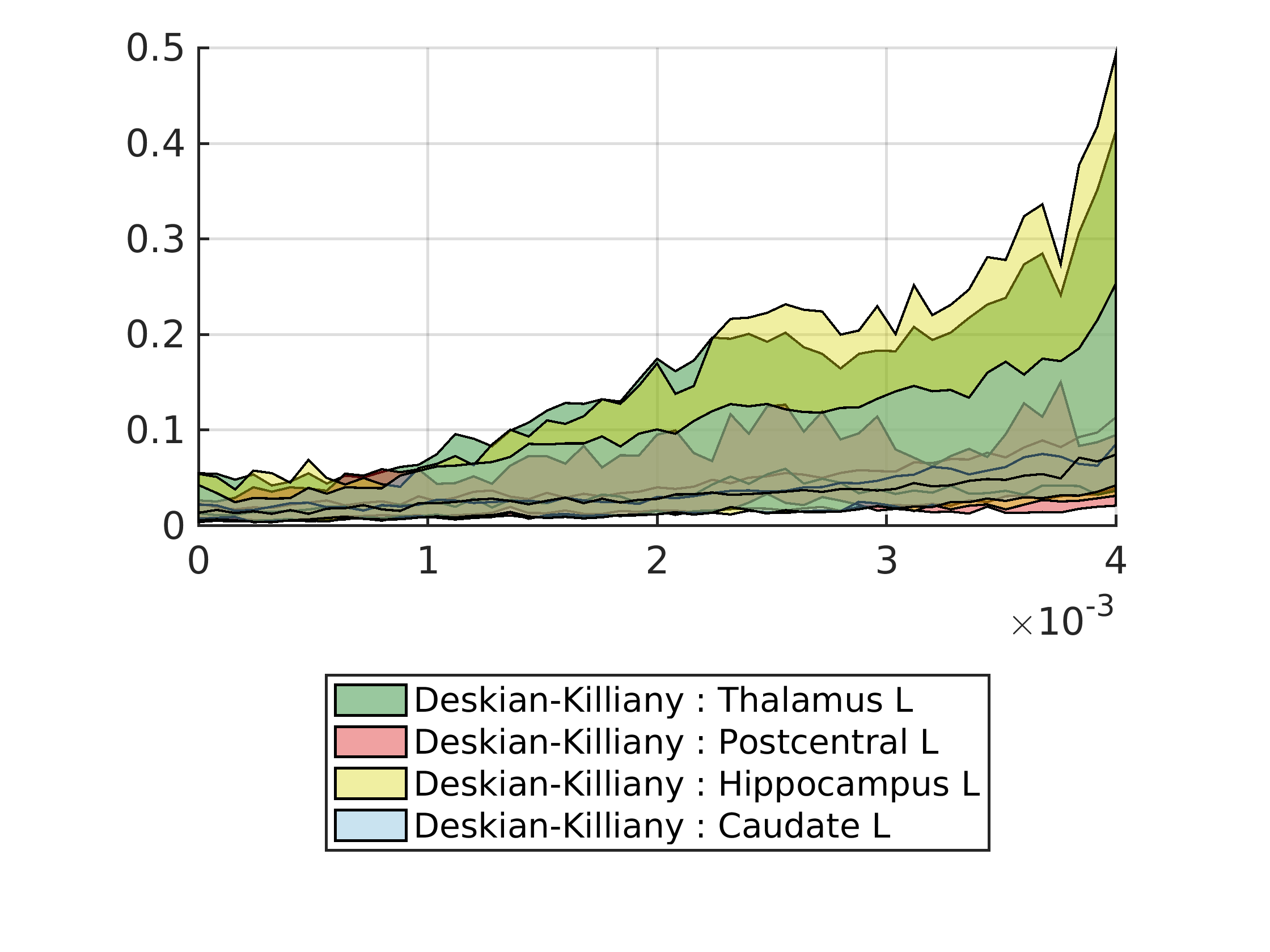}
        \caption*{EP 20, PM 20}
    \end{subfigure}
\end{minipage}

\caption{Comparison of the effect of prior parameters EP-SNR and PM-SNR on the source reconstruction using time series plots for non-Smoothed (left) and Smoothed with standardization exponent 1.25 (right) for noise level 10 db. The colors are as in Figure \ref{fig:source_placement}. }
\label{fig:comparison_noise10}
\end{figure*}

\begin{figure*}[h!]
\centering
\noindent
\begin{minipage}[t]{0.46\textwidth}
    \centering
    \begin{subfigure}{0.48\textwidth}
        \includegraphics[width=\linewidth, trim=0 120 0 0, clip]{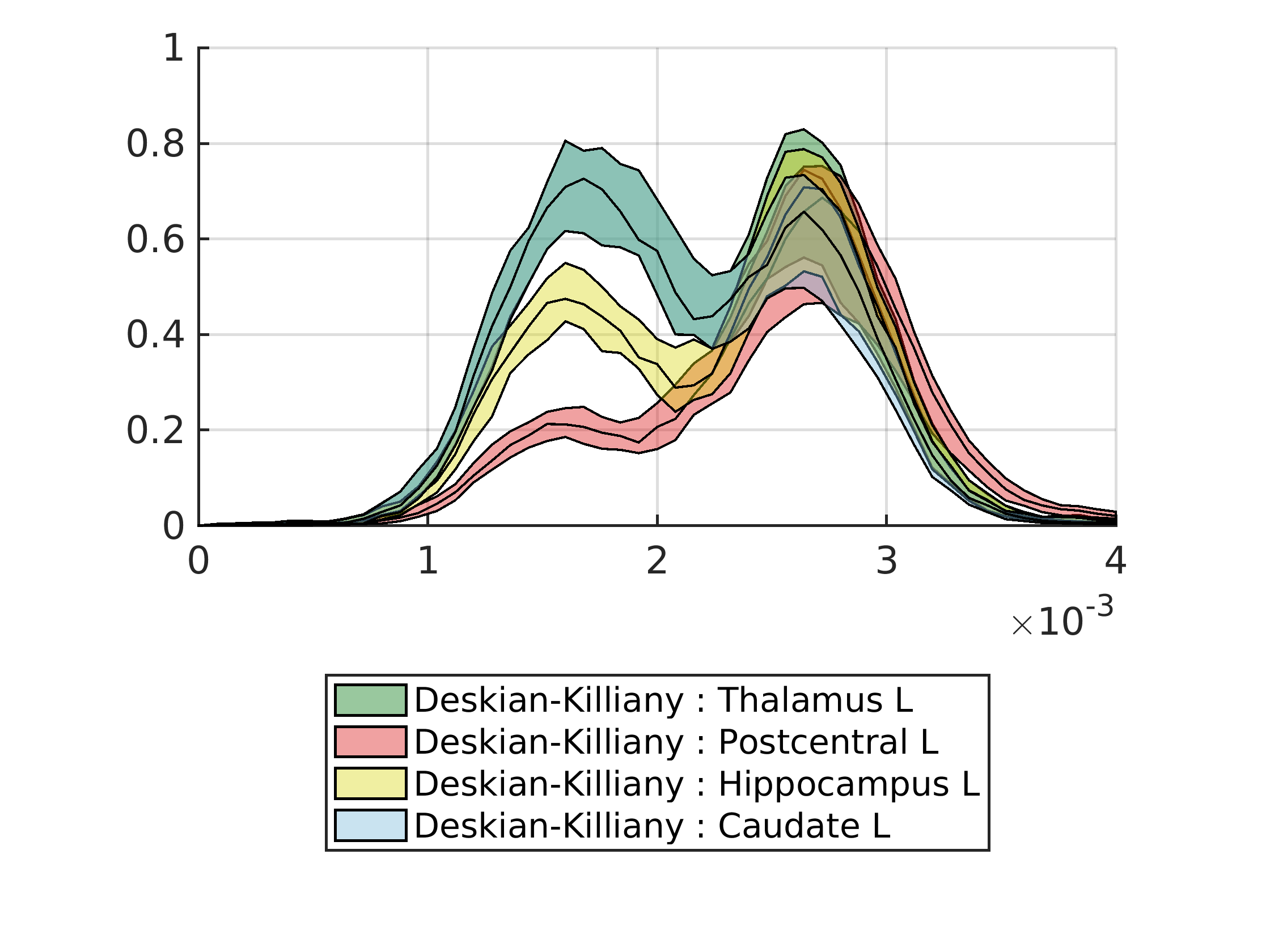}
        \caption*{EP 0, PM 0}
    \end{subfigure}
    \hfill
    \begin{subfigure}{0.48\textwidth}
        \includegraphics[width=\linewidth, trim=0 120 0 0, clip]{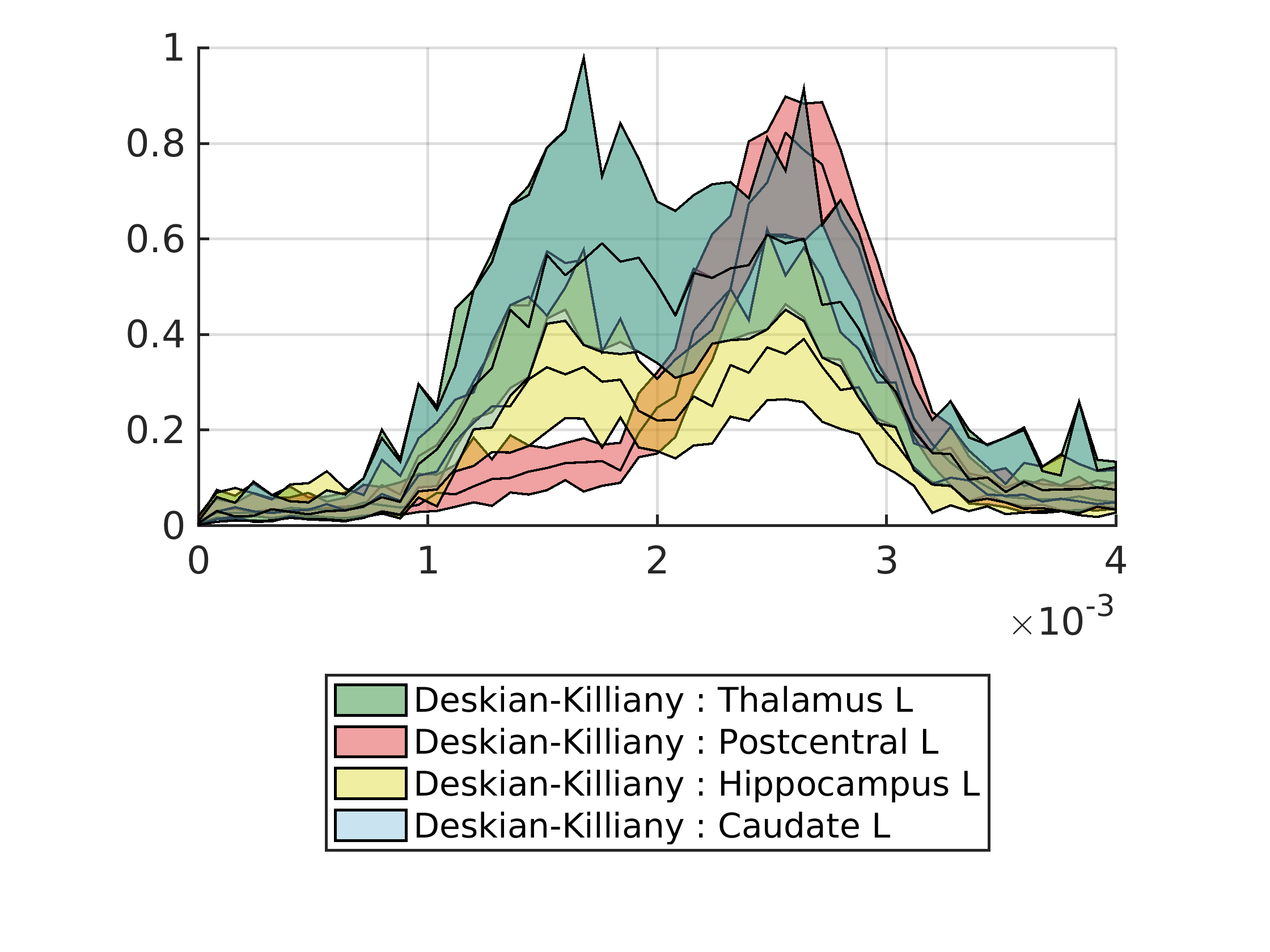}
        \caption*{EP 20, PM 0}
    \end{subfigure}

    \vspace{0.5em}

    \begin{subfigure}{0.48\textwidth}
        \includegraphics[width=\linewidth, trim=0 120 0 0, clip]{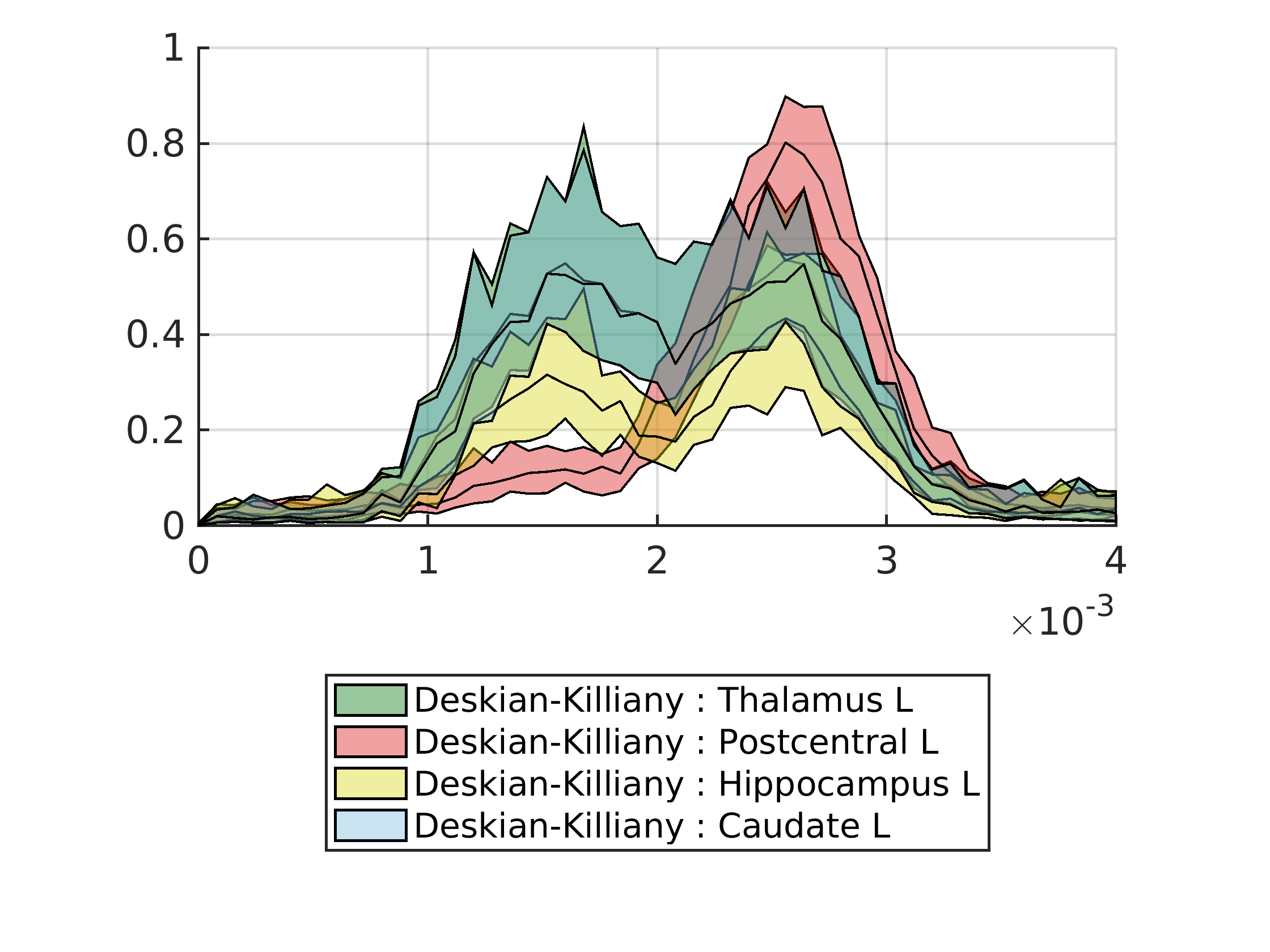}
        \caption*{EP 0, PM 20}
    \end{subfigure}
    \hfill
    \begin{subfigure}{0.48\textwidth}
        \includegraphics[width=\linewidth, trim=0 120 0 0, clip]{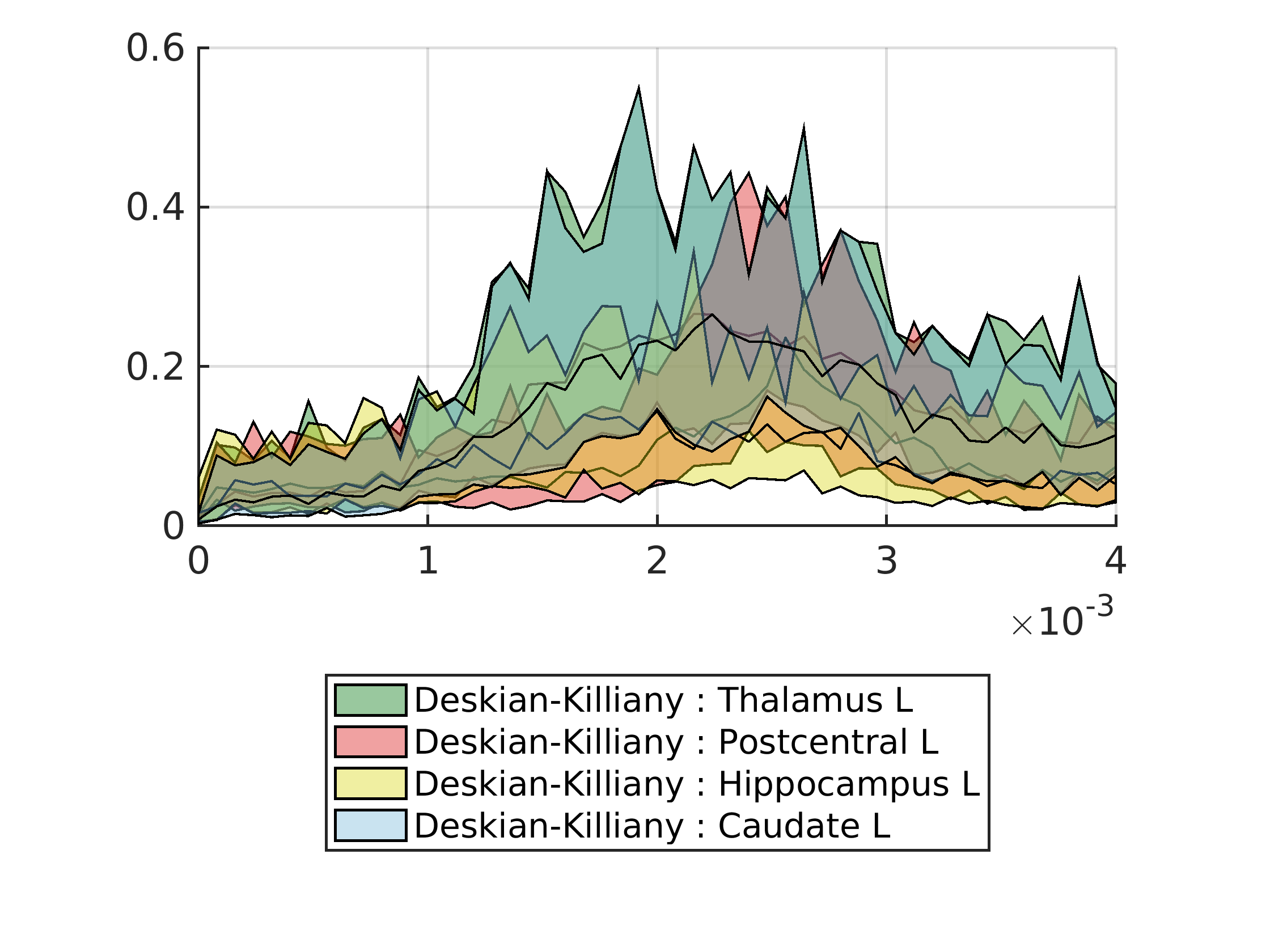}
        \caption*{EP 20, PM 20}
    \end{subfigure}
\end{minipage}
\hfill
\begin{minipage}[t]{0.46\textwidth}
    \centering
    \begin{subfigure}{0.48\textwidth}
        \includegraphics[width=\linewidth, trim=0 120 0 0, clip]{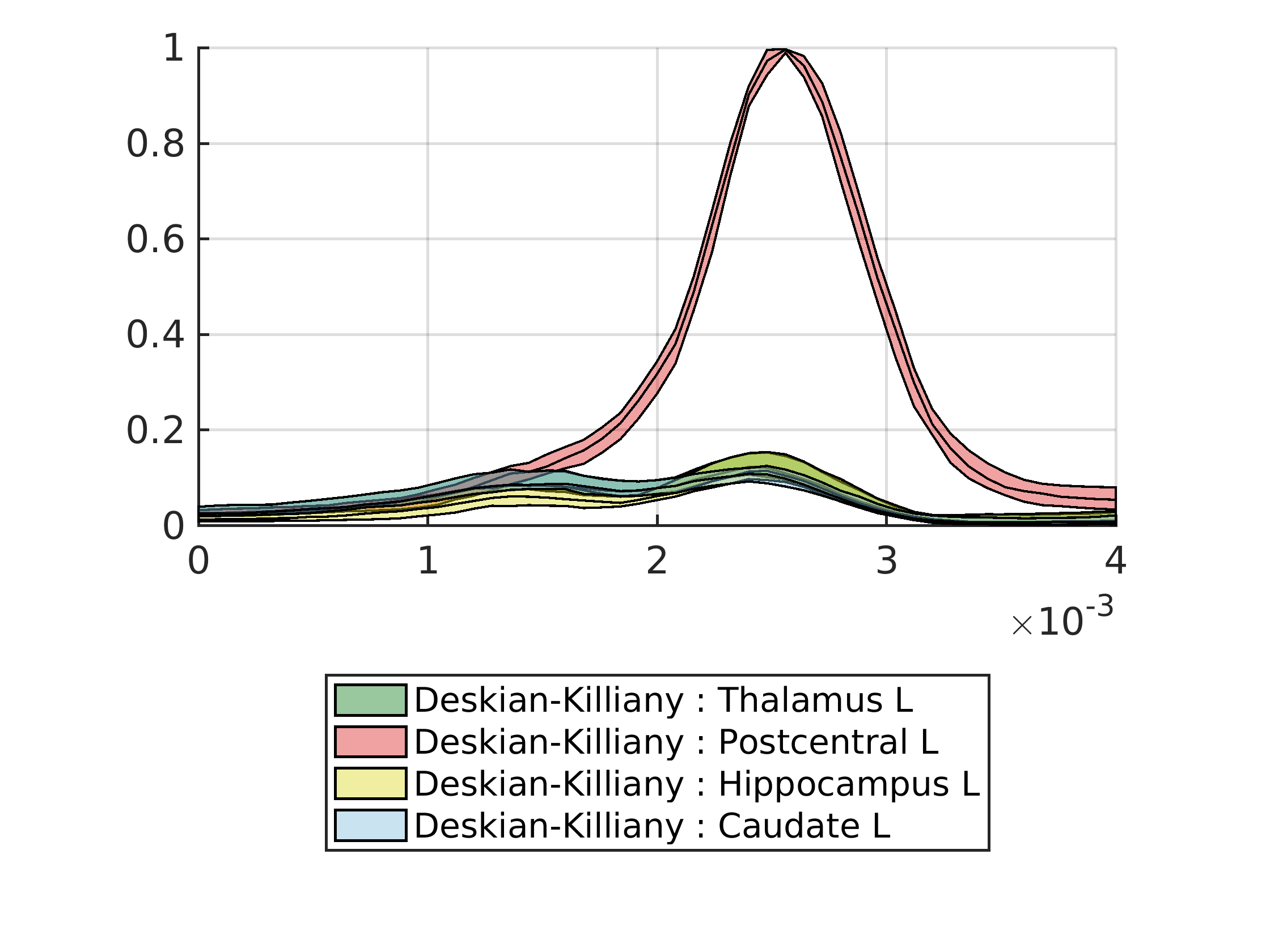}
        \caption*{EP 0, PM 0}
    \end{subfigure}
    \hfill
    \begin{subfigure}{0.48\textwidth}
        \includegraphics[width=\linewidth, trim=0 120 0 0, clip]{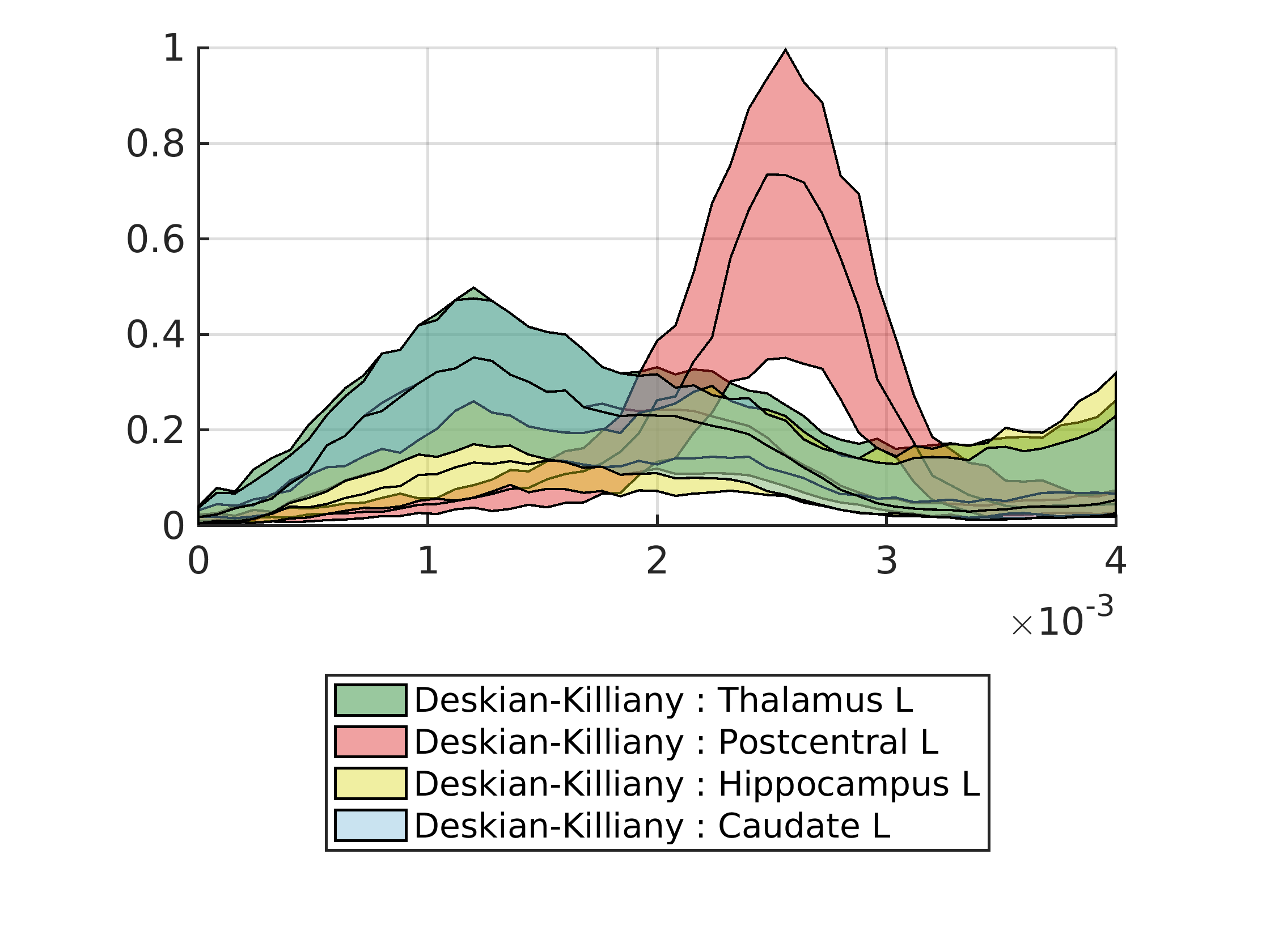}
        \caption*{EP 20, PM 0}
    \end{subfigure}

    \vspace{0.5em}

    \begin{subfigure}{0.48\textwidth}
        \includegraphics[width=\linewidth, trim=0 120 0 0, clip]{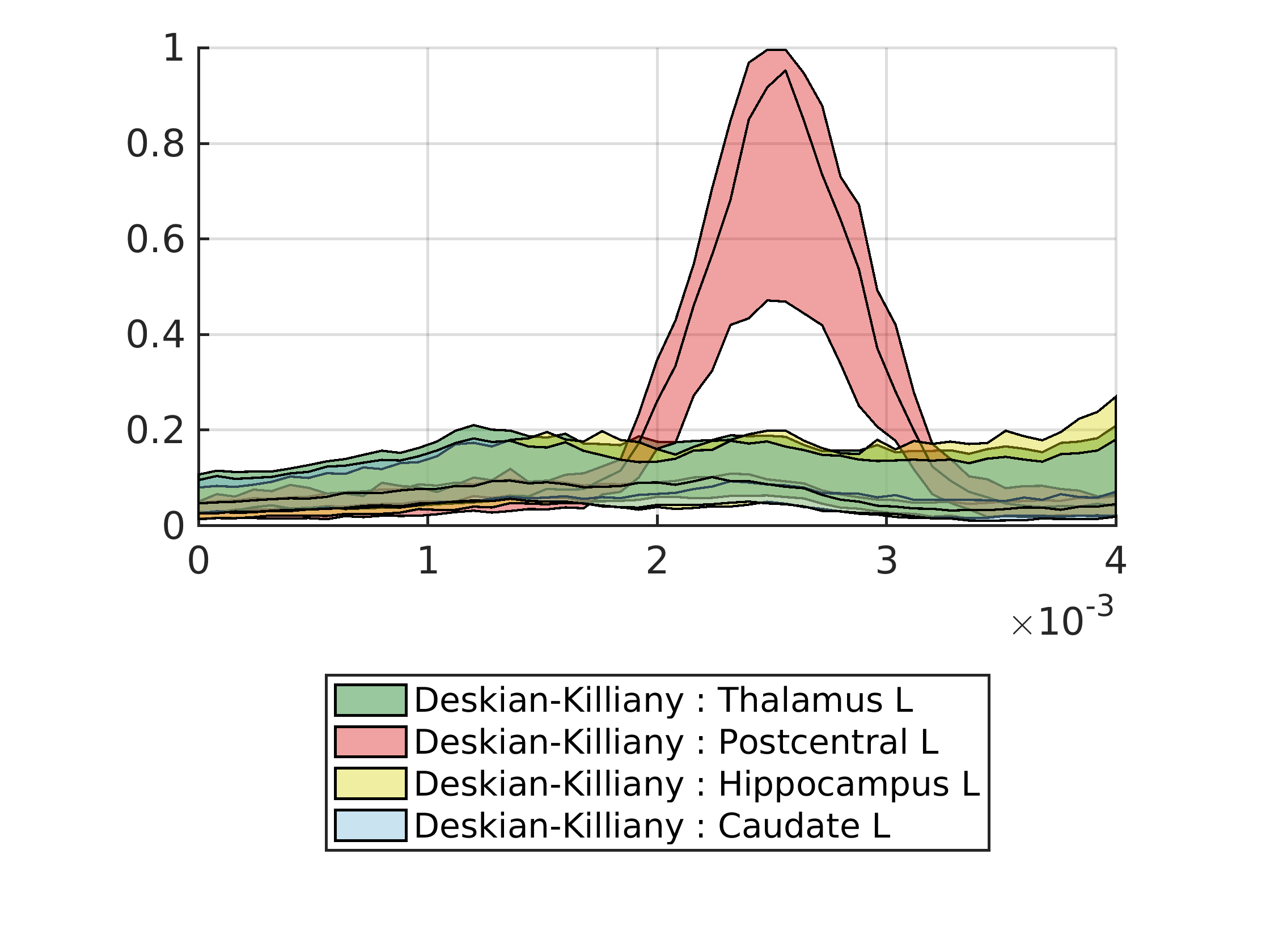}
        \caption*{EP 0, PM 20}
    \end{subfigure}
    \hfill
    \begin{subfigure}{0.48\textwidth}
        \includegraphics[width=\linewidth, trim=0 120 0 0, clip]{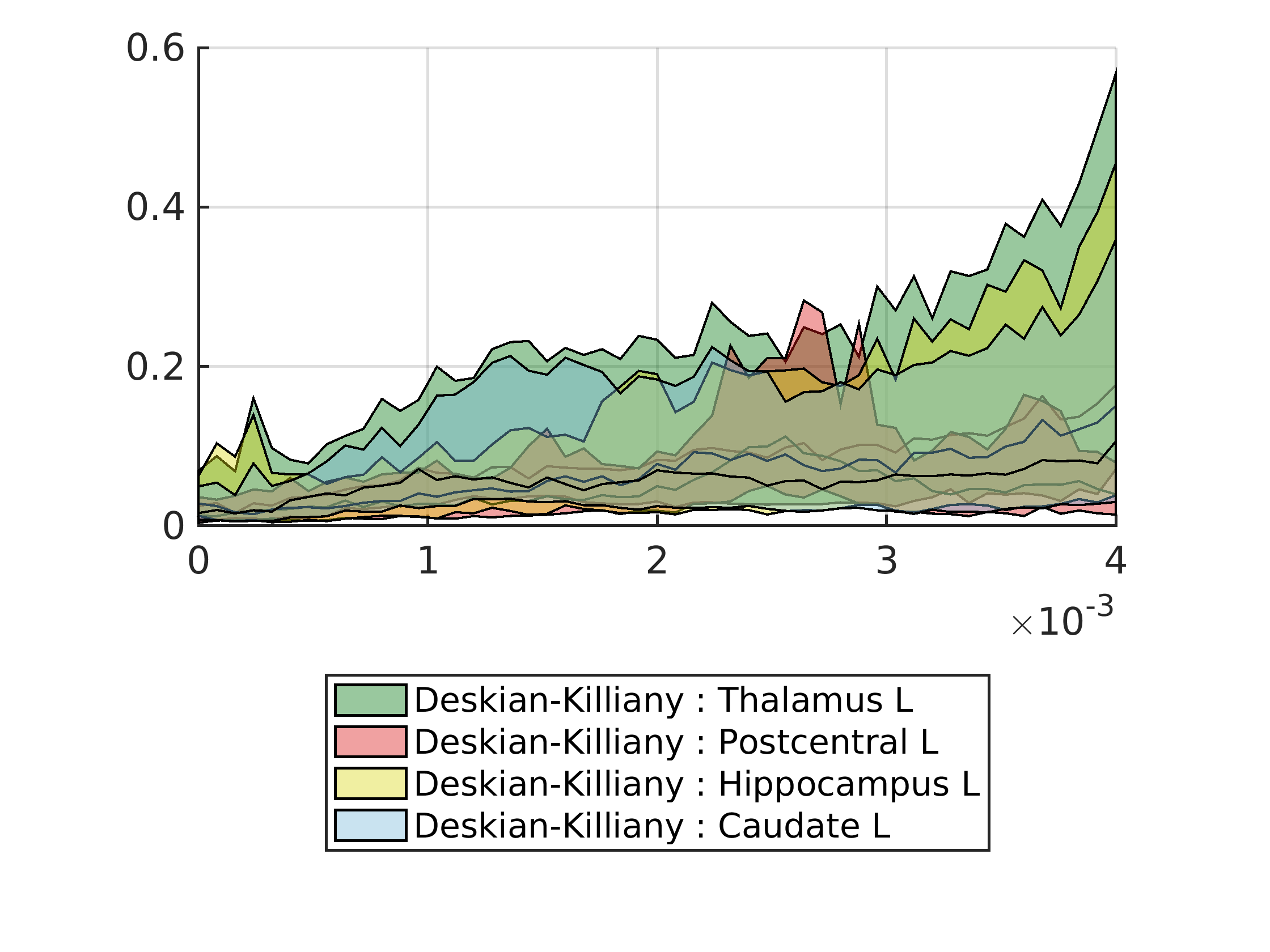}
        \caption*{EP 20, PM 20}
    \end{subfigure}
\end{minipage}
\caption{Comparison of the effect of prior parameters EP-SNR and PM-SNR on the source reconstruction using time series plots for non-Smoothed (left) and Smoothed with standardization exponent 1.25 (right) for noise level 20 dB. The colors are as in Figure \ref{fig:source_placement}.}
\label{fig:comparison_noise20}
\end{figure*}

\begin{figure*}[h!]
\centering
\noindent
\begin{minipage}[t]{0.46\textwidth}
    \centering
    \begin{subfigure}{0.48\textwidth}
        \includegraphics[width=\linewidth, trim=0 120 0 0, clip]{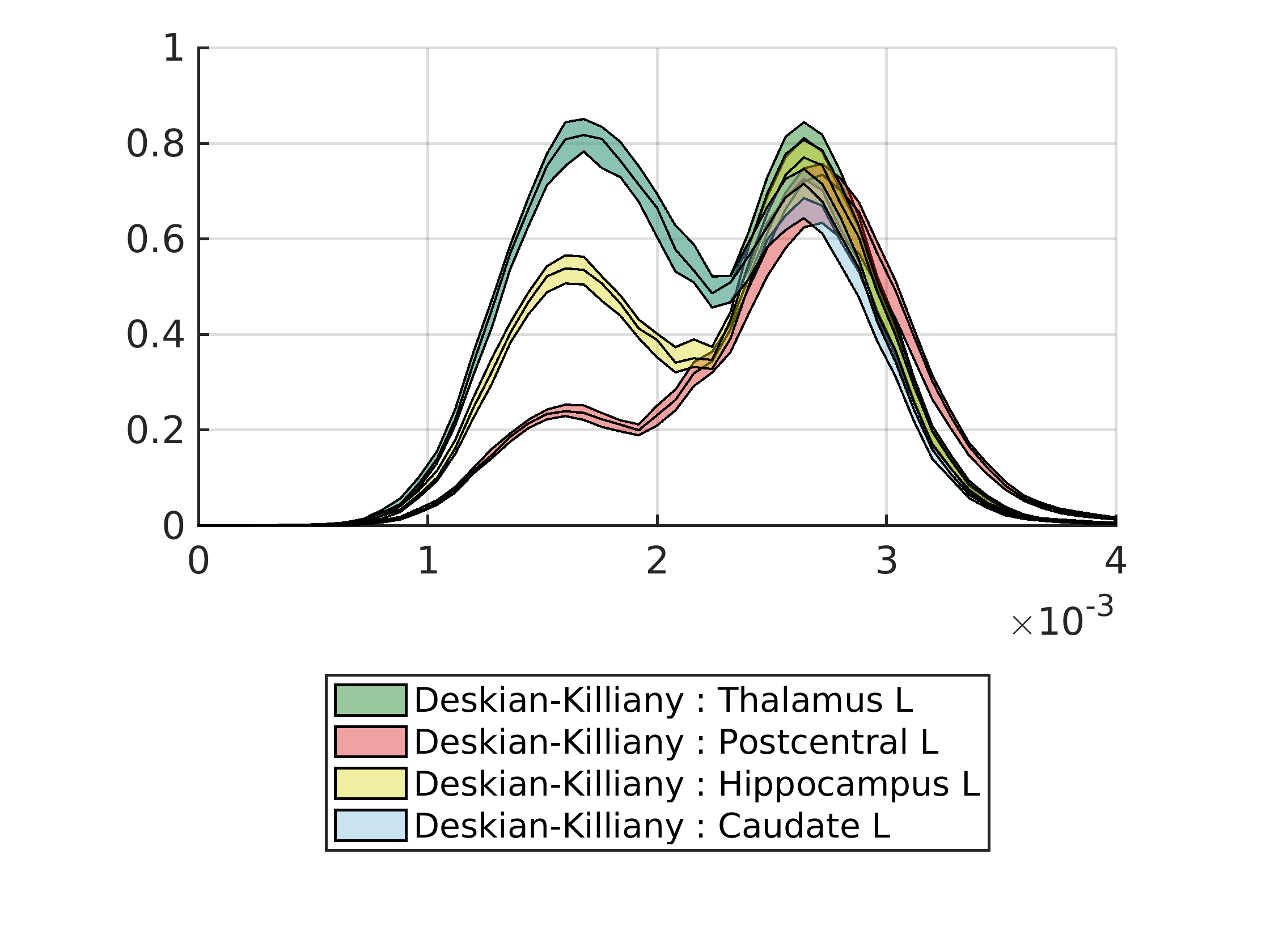}
        \caption*{EP 0, PM 0}
    \end{subfigure}
    \hfill
    \begin{subfigure}{0.48\textwidth}
        \includegraphics[width=\linewidth, trim=0 120 0 0, clip]{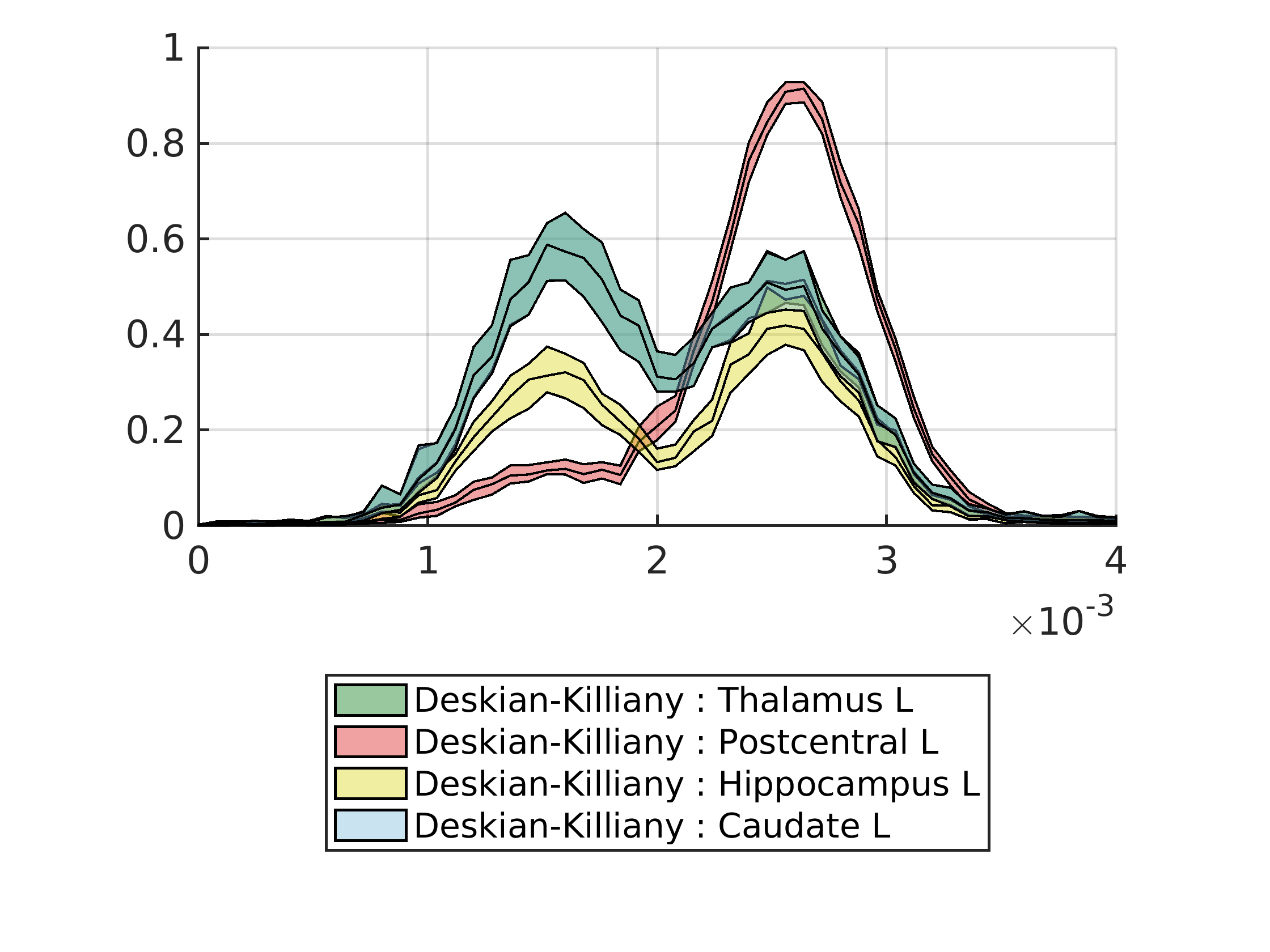}
        \caption*{EP 20, PM 0}
    \end{subfigure}

    \vspace{0.5em}

    \begin{subfigure}{0.48\textwidth}
        \includegraphics[width=\linewidth, trim=0 120 0 0, clip]{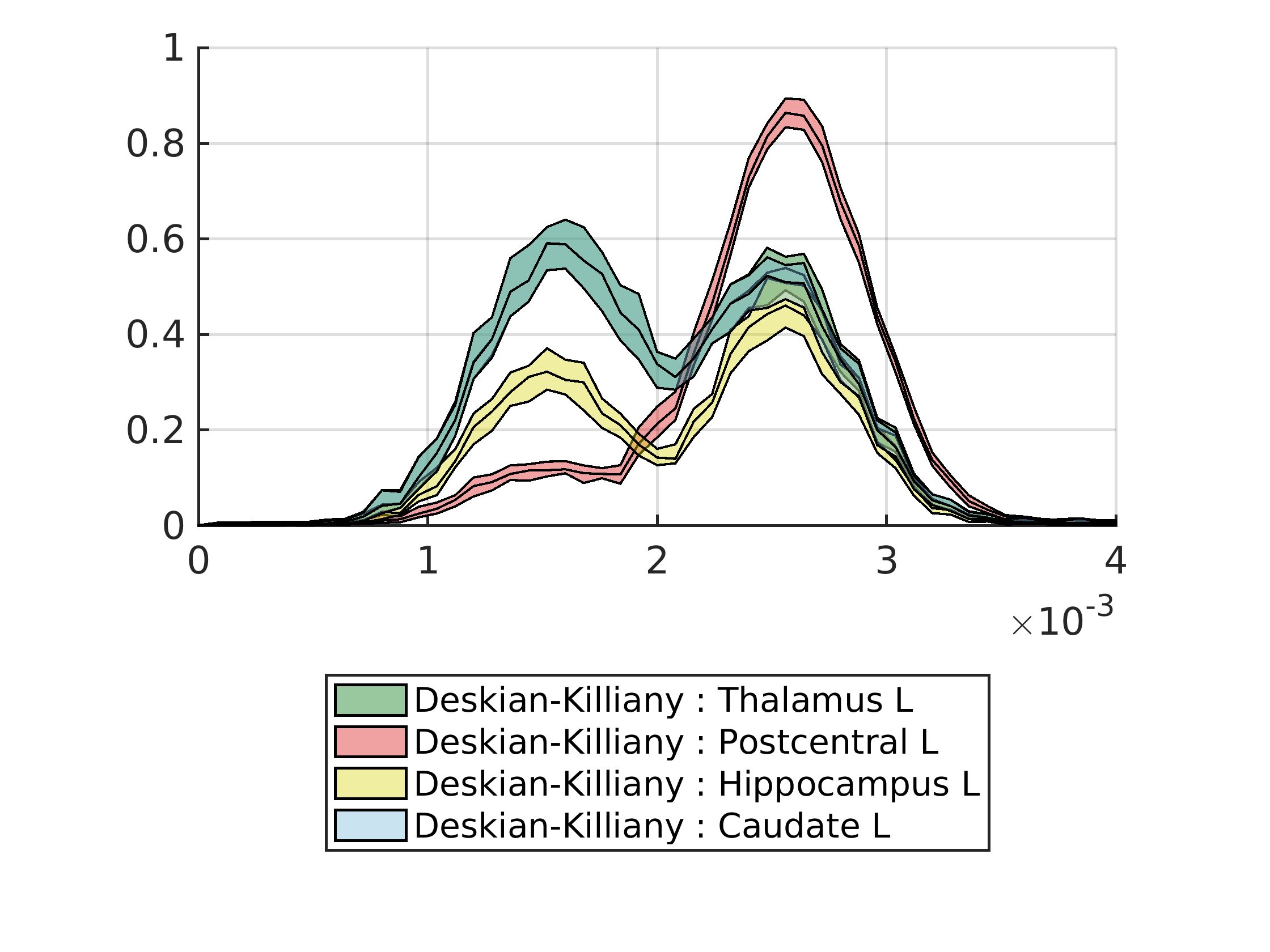}
        \caption*{EP 0, PM 20}
    \end{subfigure}
    \hfill
    \begin{subfigure}{0.48\textwidth}
        \includegraphics[width=\linewidth, trim=0 120 0 0, clip]{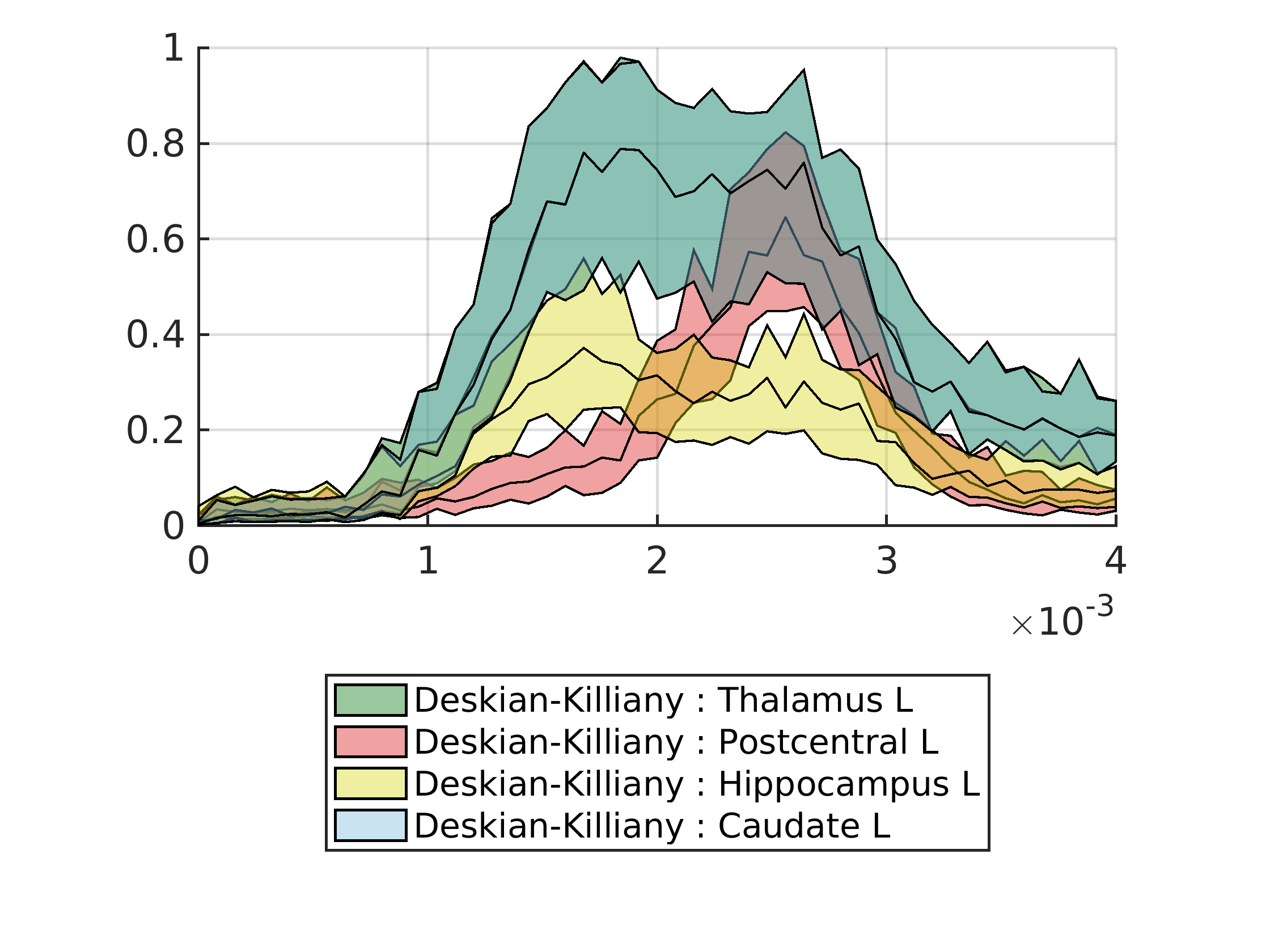}
        \caption*{EP 20, PM 20}
    \end{subfigure}
\end{minipage}
\hfill
\begin{minipage}[t]{0.46\textwidth}
    \centering
    \begin{subfigure}{0.48\textwidth}
        \includegraphics[width=\linewidth, trim=0 120 0 0, clip]{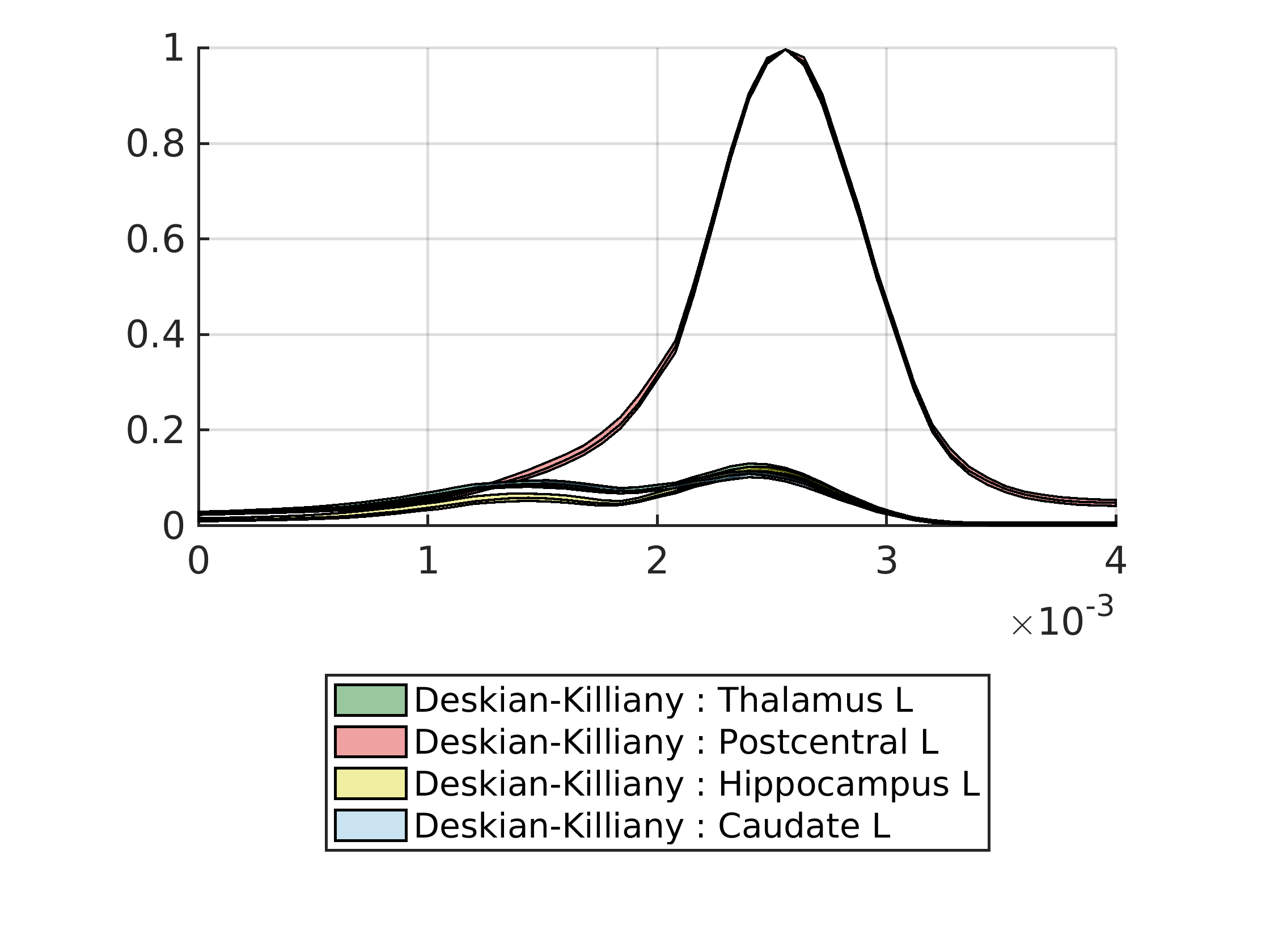}
        \caption*{EP 0, PM 0}
    \end{subfigure}
    \hfill
    \begin{subfigure}{0.48\textwidth}
        \includegraphics[width=\linewidth, trim=0 120 0 0, clip]{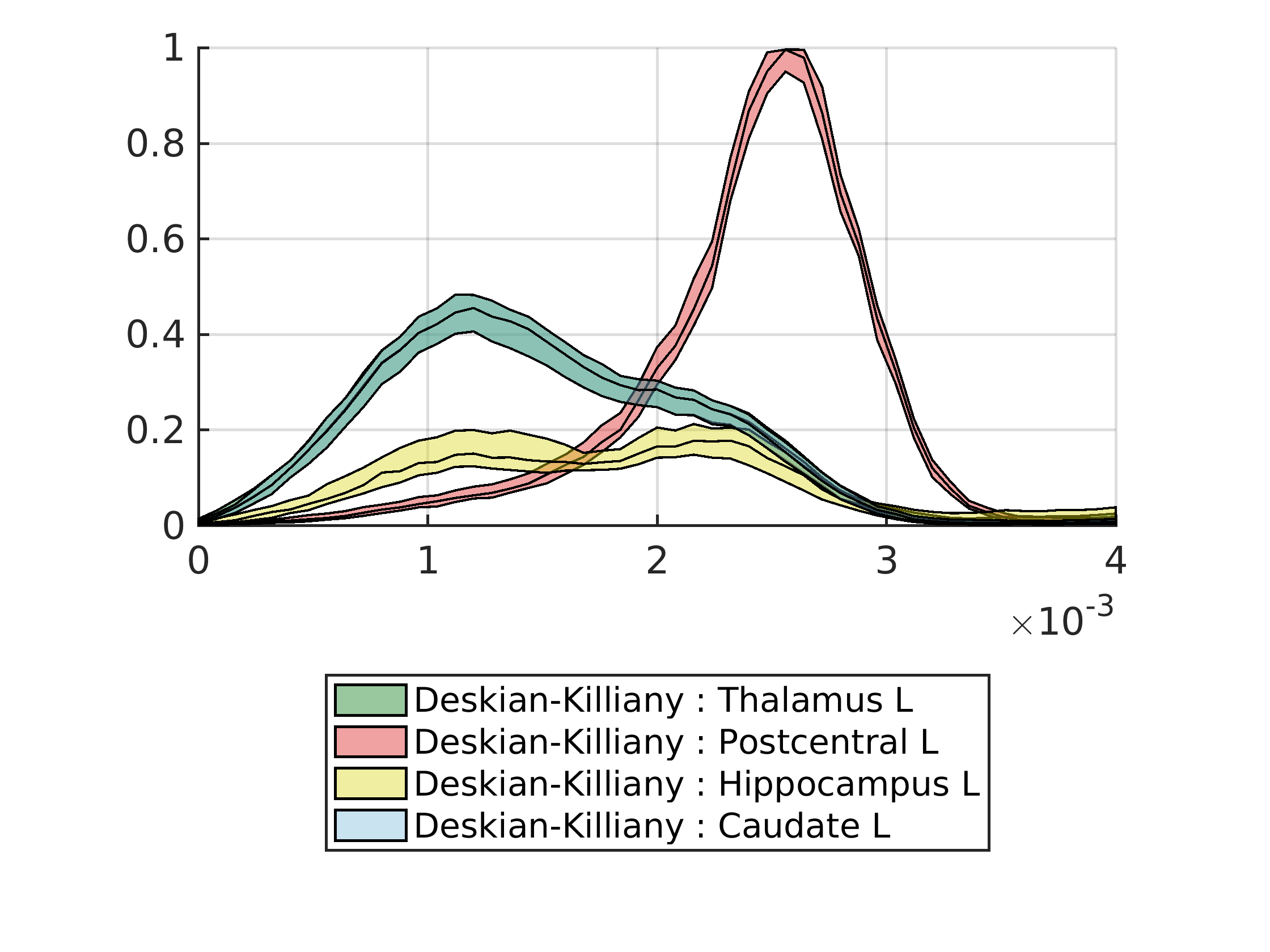}
        \caption*{EP 20, PM 0}
    \end{subfigure}

    \vspace{0.5em}

    \begin{subfigure}{0.48\textwidth}
        \includegraphics[width=\linewidth, trim=0 120 0 0, clip]{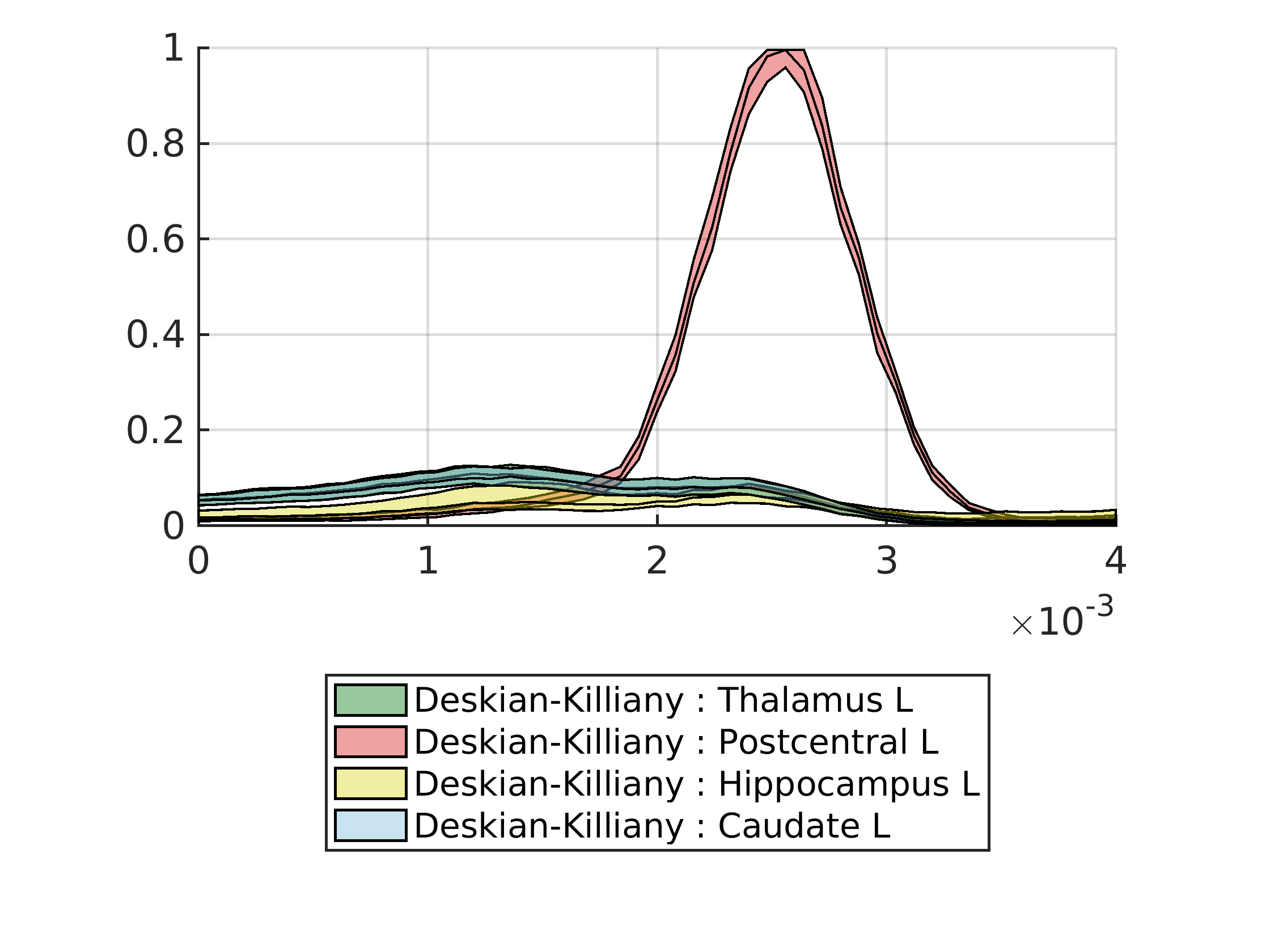}
        \caption*{EP 0, PM 20}
    \end{subfigure}
    \hfill
    \begin{subfigure}{0.48\textwidth}
        \includegraphics[width=\linewidth, trim=0 120 0 0, clip]{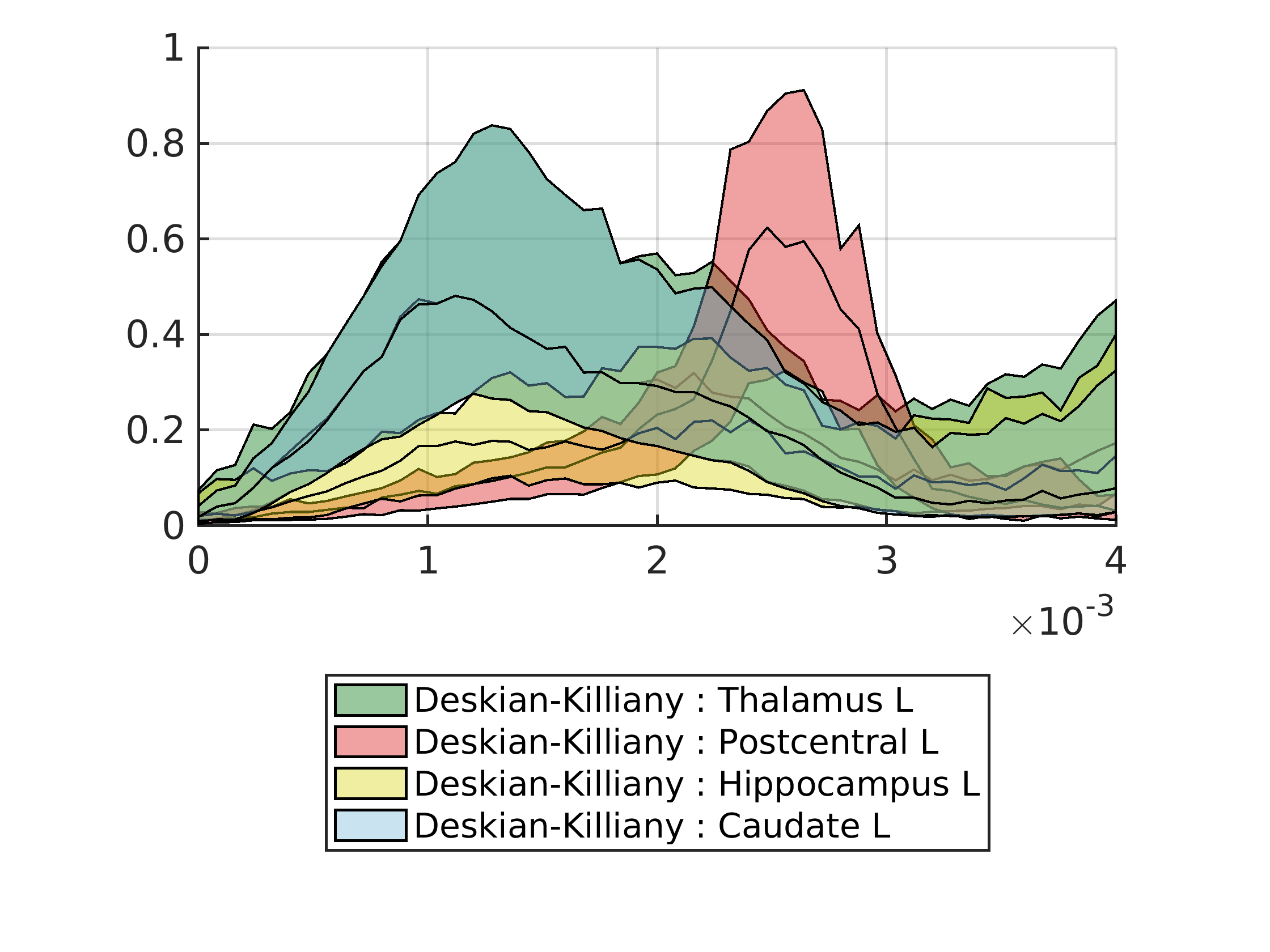}
        \caption*{EP 20, PM 20}
    \end{subfigure}
\end{minipage}
\caption{Comparison of the effect of prior parameters EP-SNR and PM-SNR on the source reconstruction using time series plots for non-Smoothed (left) and Smoothed with standardization exponent 1.25 (right) for noise level 30 dB. The colors are as in Figure \ref{fig:source_placement}.}
\label{fig:comparison_noise30}
\end{figure*}

The effect of smoothing was assessed at 10 dB, 20 dB, and 30 dB noise levels by comparing time series plots between the smoothed and non-smoothed cases, using the standardization exponent 1.25. At each noise level, EP-SNR 20 and PM-SNR 0 were found to produce the best results when combined with smoothing. Smoothing caused a visible separation of the deep and superficial activities, although at 10 dB and 20 dB, the smoothed time series appeared slightly more volatile. At 30 dB, the effect of smoothing was clearly beneficial, leading to visibly cleaner time courses (Figure \ref{fig:comparison_noise10}, Figure \ref{fig:comparison_noise20} and Figure \ref{fig:comparison_noise30}).

\begin{figure*}[ht!]
\centering
\noindent
\begin{minipage}[t]{0.46\textwidth}
    \centering
    \begin{subfigure}{0.48\textwidth}
        \includegraphics[width=\linewidth]{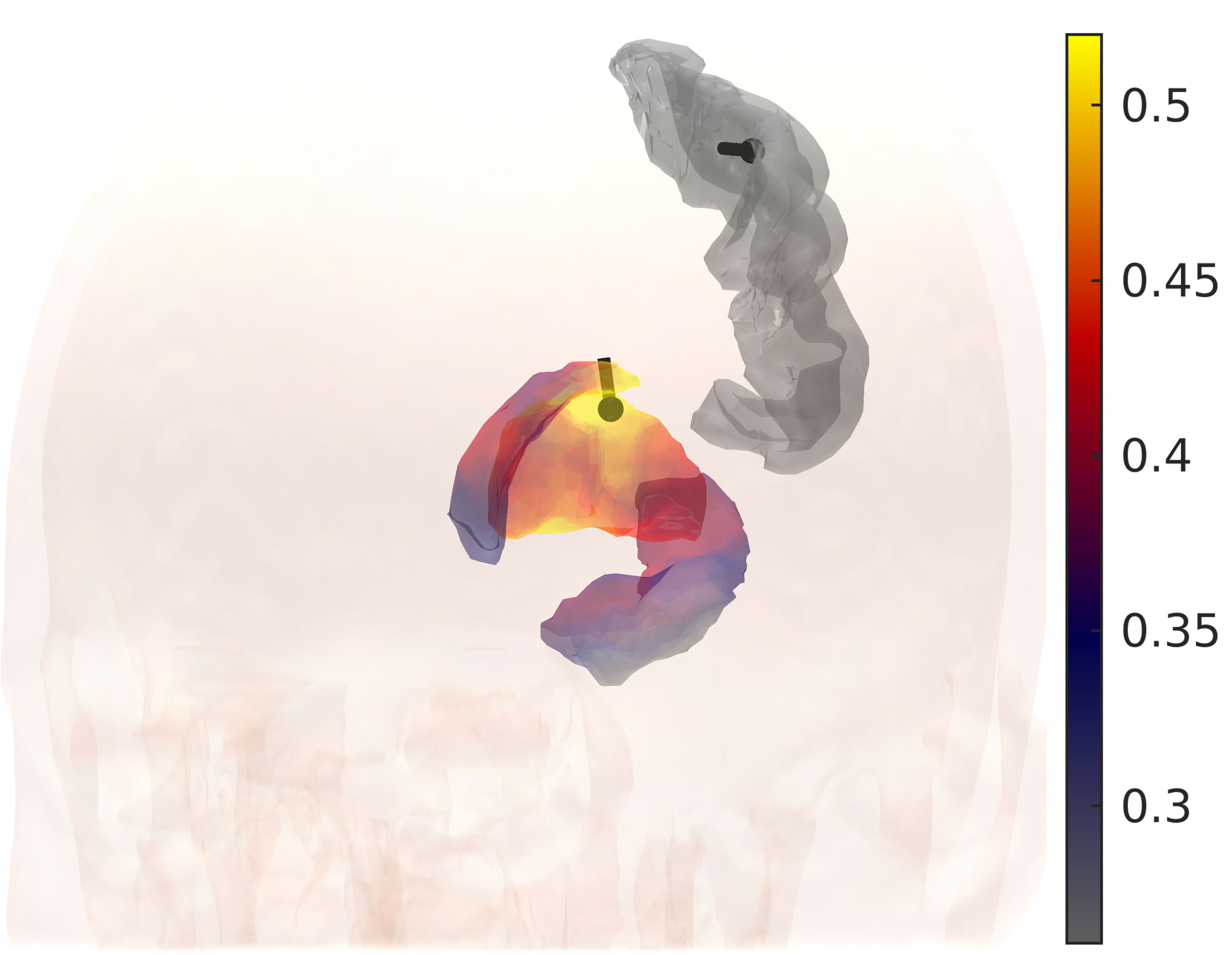}
        \caption*{10 dB 1.5 ms}
    \end{subfigure}
    \hfill
    \begin{subfigure}{0.48\textwidth}
        \includegraphics[width=\linewidth]{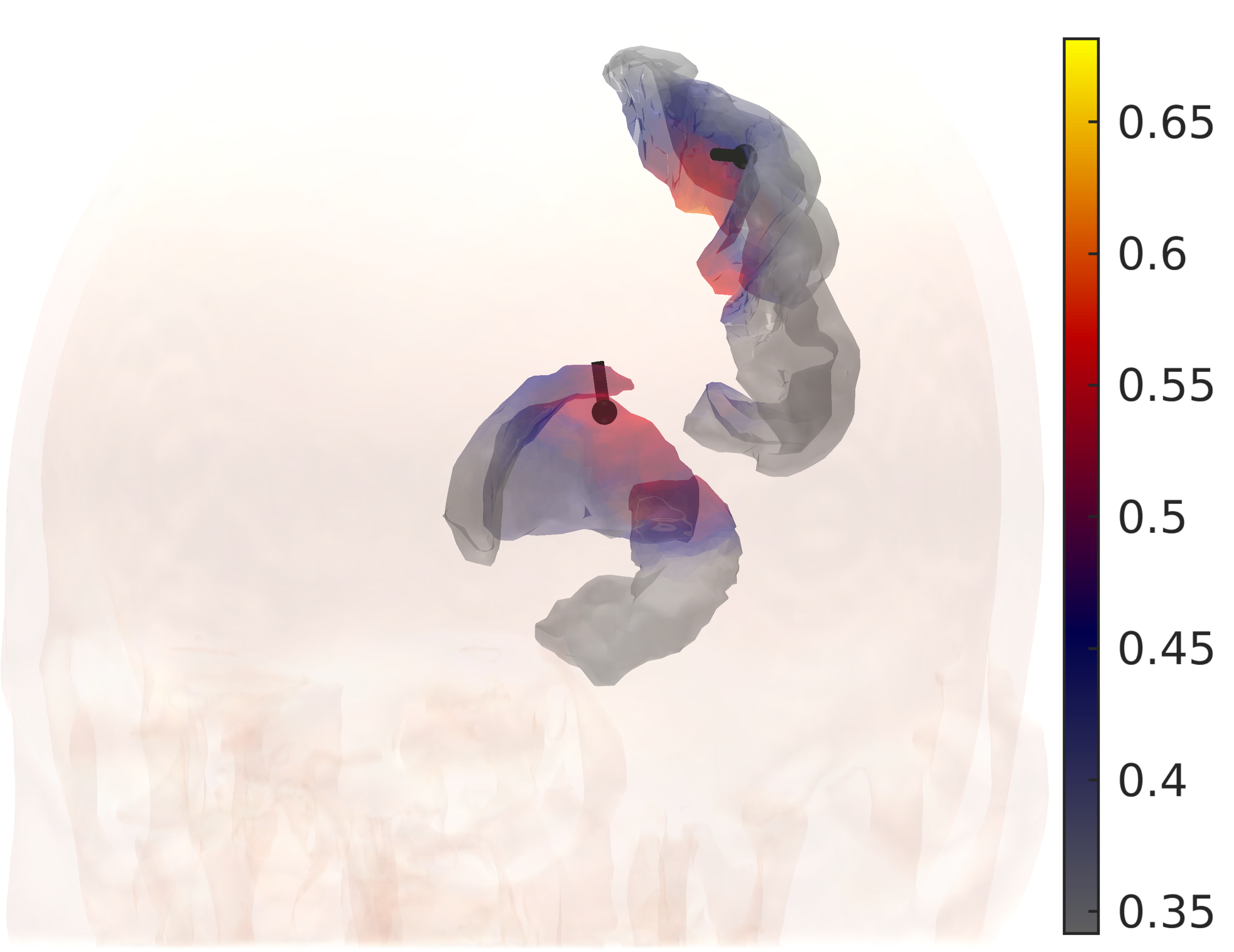}
        \caption*{10 dB 2.5 ms}
    \end{subfigure}

    \vspace{0.5em}

    \begin{subfigure}{0.48\textwidth}
        \includegraphics[width=\linewidth]{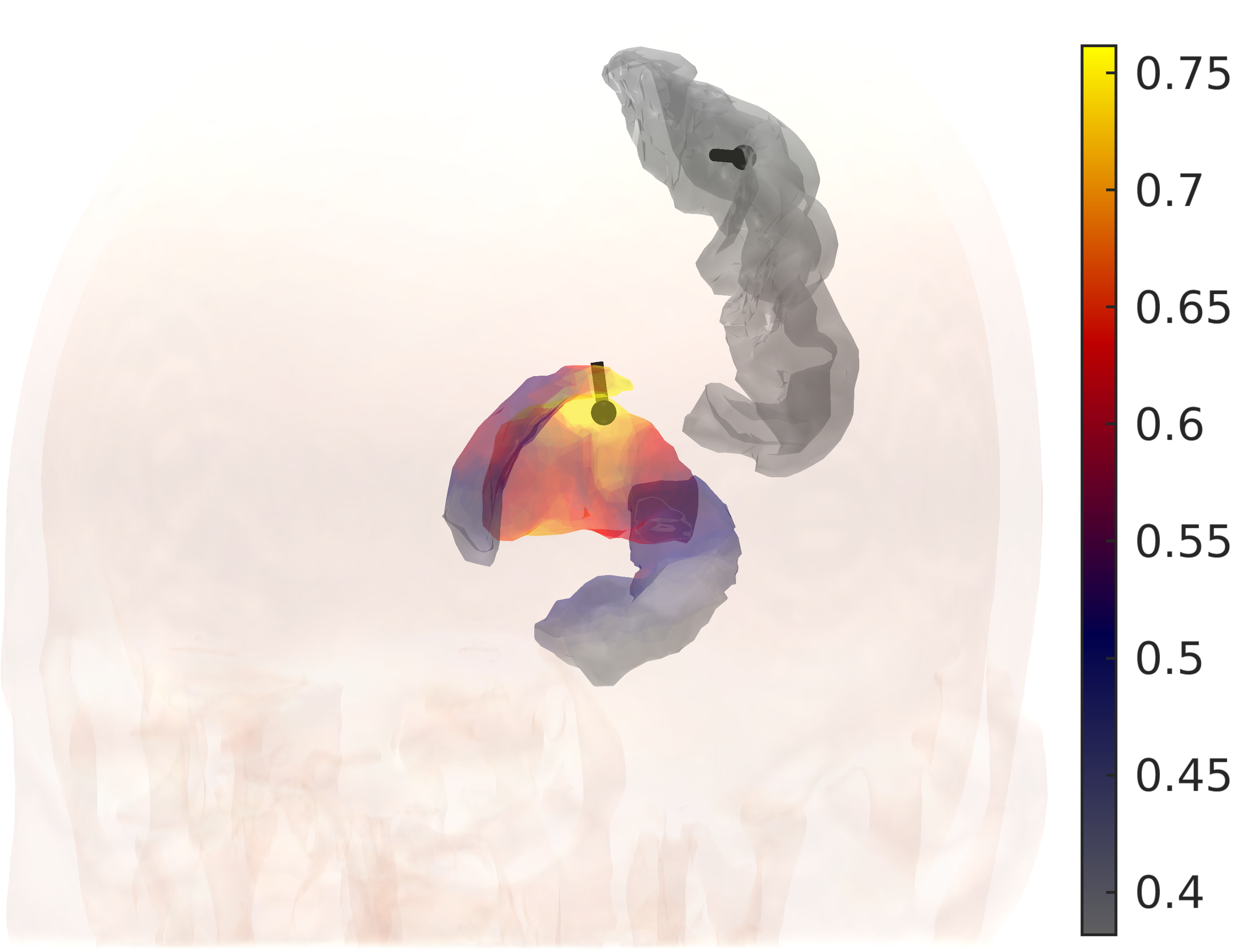}
        \caption*{30 dB 1.5 ms}
    \end{subfigure}
    \hfill
    \begin{subfigure}{0.48\textwidth}
        \includegraphics[width=\linewidth]{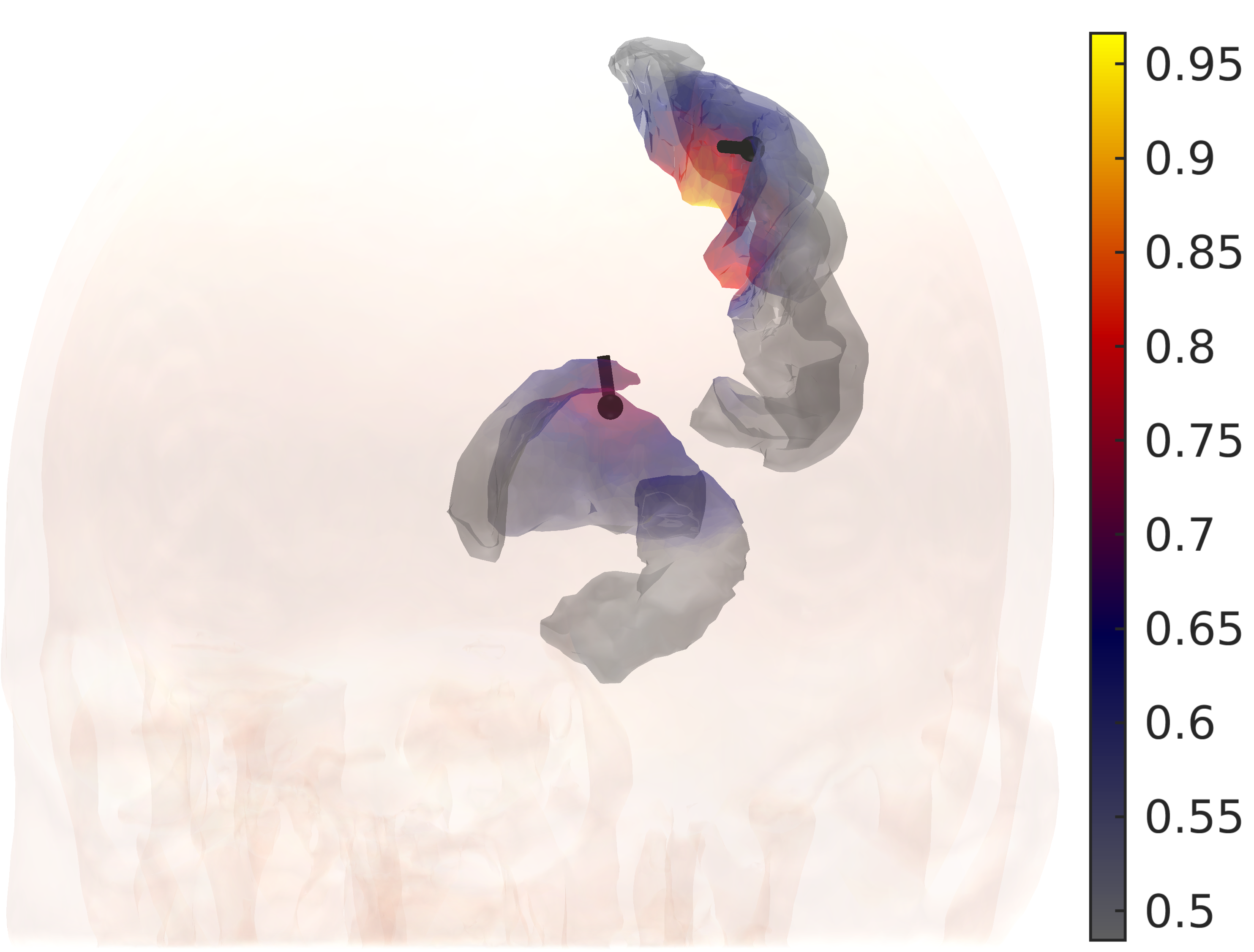}
        \caption*{30 dB 2.5 ms}
    \end{subfigure}
\end{minipage}
\hfill
\begin{minipage}[t]{0.46\textwidth}
    \centering
    \begin{subfigure}{0.48\textwidth}
        \includegraphics[width=\linewidth]{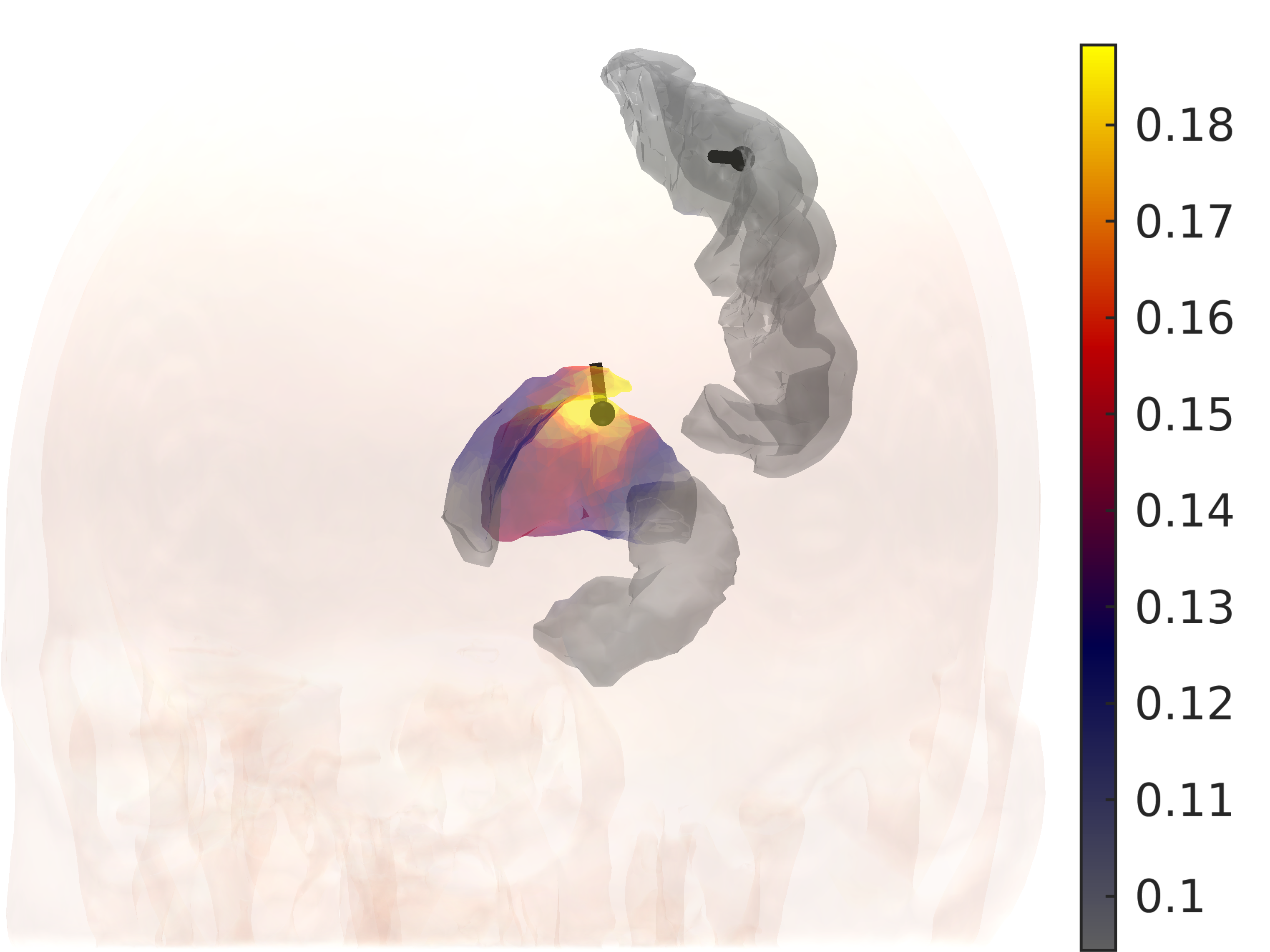}
        \caption*{10 dB 1.5 ms}
    \end{subfigure}
    \hfill
    \begin{subfigure}{0.48\textwidth}
        \includegraphics[width=\linewidth]{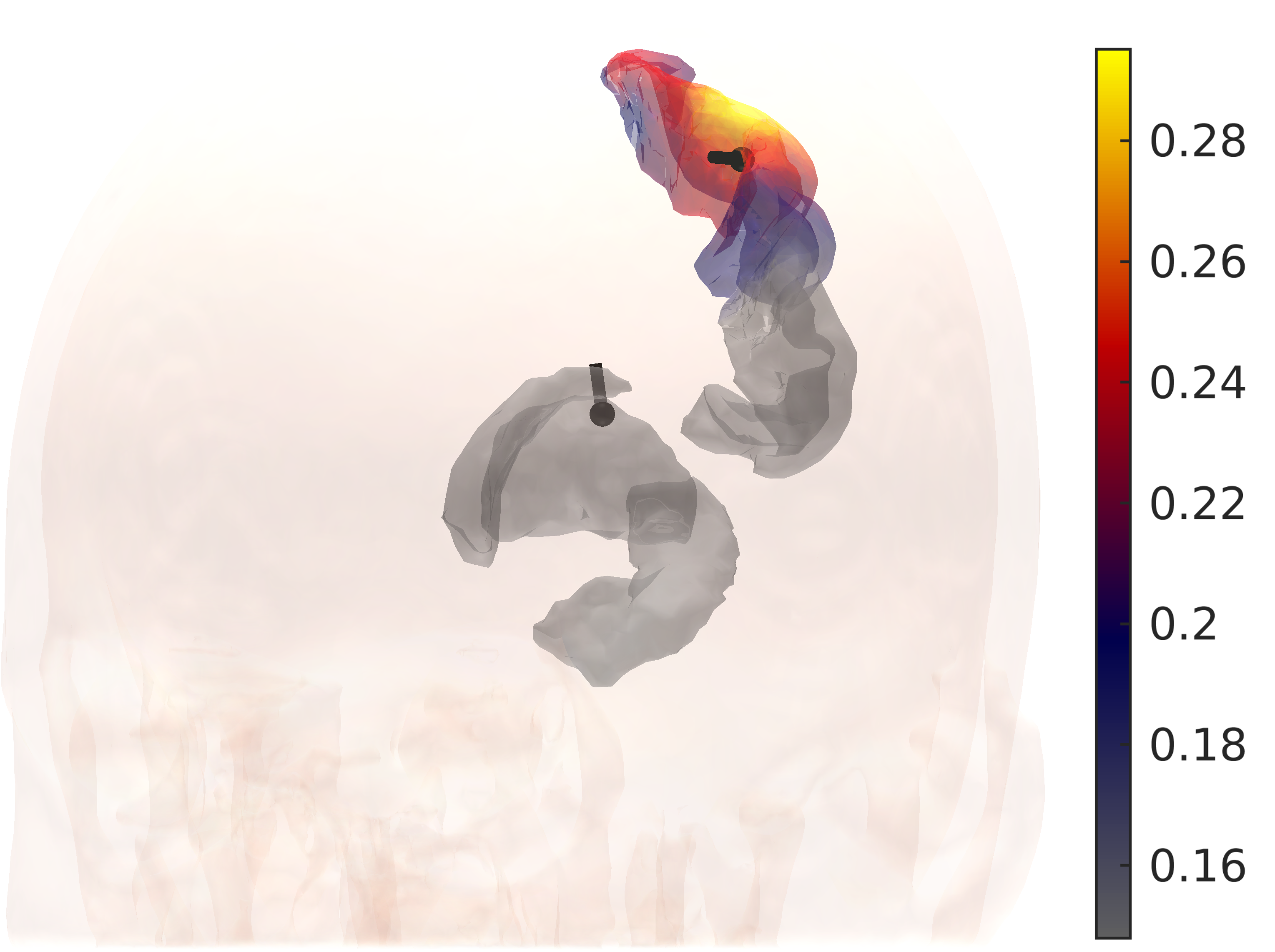}
        \caption*{10 dB 2.5 ms}
    \end{subfigure}

    \vspace{0.5em}

    \begin{subfigure}{0.48\textwidth}
        \includegraphics[width=\linewidth]{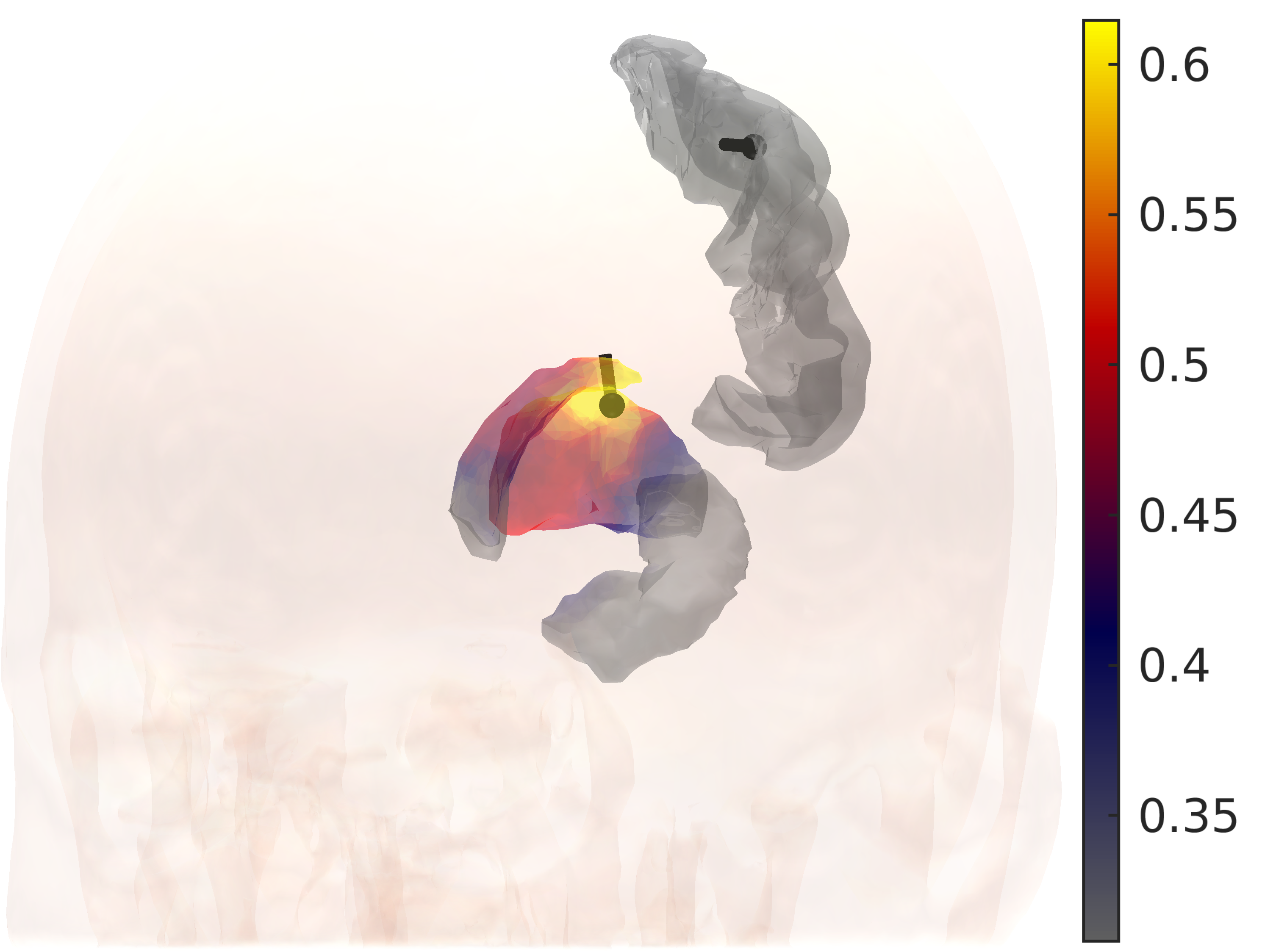}
        \caption*{30 dB 1.5 ms}
    \end{subfigure}
    \hfill
    \begin{subfigure}{0.48\textwidth}
        \includegraphics[width=\linewidth]{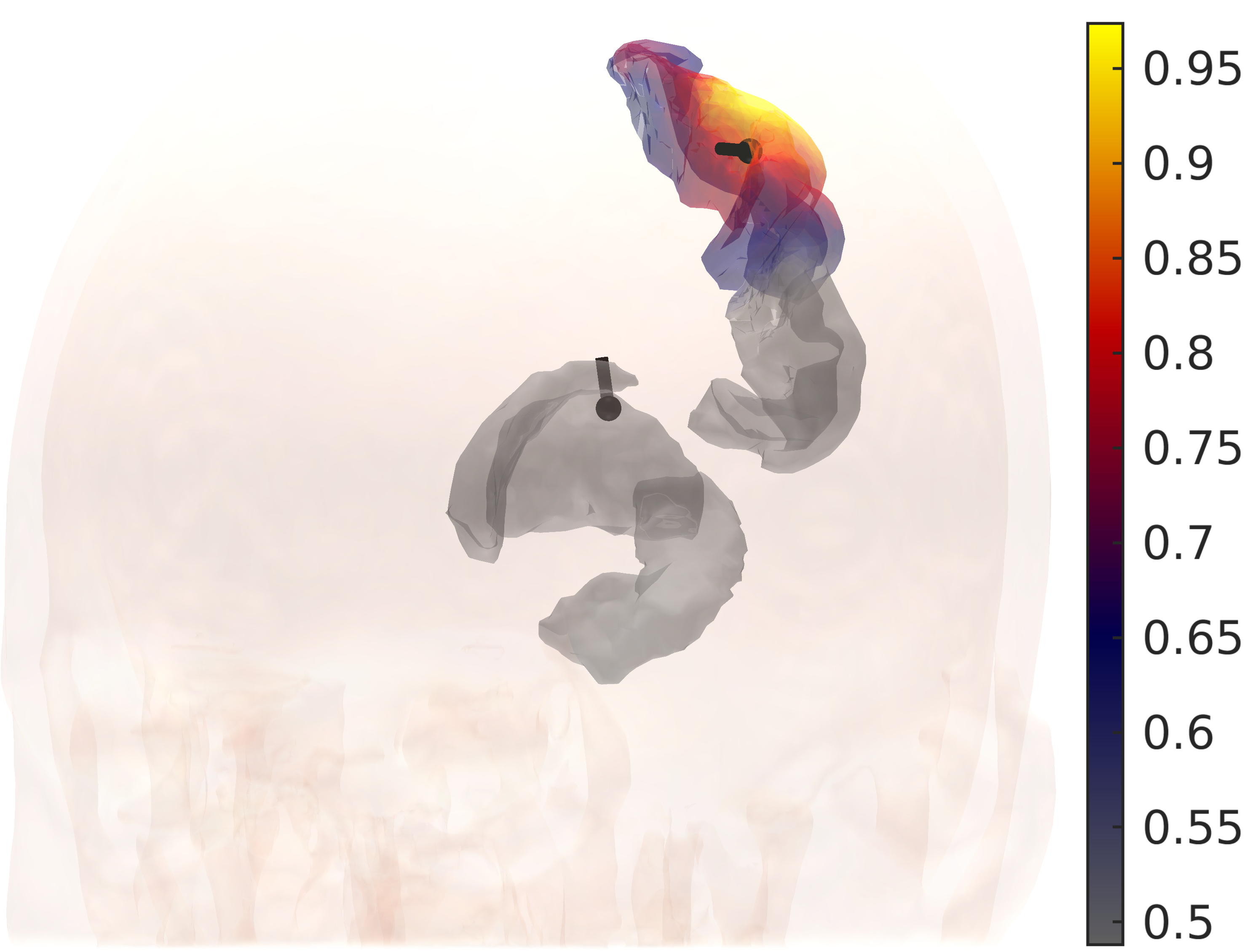}
        \caption*{30 dB 2.5 ms}
    \end{subfigure}
\end{minipage}

\caption{Comparison of average reconstructions for non-smoothed (left) and smoothed case with standardization exponent 1.25 (right) at the time of deep activation (1.5 ms) and superficial activation (2.5 ms) with high noise (10 dB) and low noise (30 dB) cases. The measurement peak-SNR is limited to 10-30 dB in these experiments. We assume that the dynamic range of a reconstruction stays on a lower level due to limiting factors, such as the lead field attenuation and modeling uncertainty; that is why we have limited the low-end of the colormap to -6 dB with respect to the maximum of the reconstruction (the dynamic range of the visualization is 6 dB).}
\label{fig:average_reconstructions}
\end{figure*}

The average reconstructions at the two critical time points revealed the most significant improvement with smoothing. In the non-smoothed reconstructions, residual deep activity could be seen at the time of the superficial peak at all noise levels, indicating incomplete separation between the sources. After smoothing along with the standardization exponent 1.25, this residual deep activity was successfully suppressed, resulting in a clear separation of the deep and superficial activations. This improvement was consistent across all noise levels considered, demonstrating the advantage of combining RTS smoothing with standardized Kalman filtering along with a suitable standardization exponent (Figure \ref{fig:average_reconstructions}).

\section{Discussion}

The main findings of this study demonstrate that careful tuning of prior parameters in the SKF framework significantly enhances EEG source reconstruction accuracy, particularly for deep brain activity. Specifically, we found that optimizing EP-SNR and PM-SNR leads to more accurate estimation of both deep and superficial sources. Moreover, incorporating a tuned standardization exponent in conjunction with RTS smoothing produced the most robust and spatially precise reconstructions. This combination improved the separability of overlapping sources and reduced depth bias, especially under low-noise conditions. These findings validate the utility of SKF in modeling temporally dynamic and spatially extended brain activity, supporting its application in scenarios where deep source detection is critical.

The findings of this study highlight the potential of the SKF as a dynamic and flexible tool for EEG source localization, particularly for reconstructing deep brain activity. By demonstrating that careful tuning of prior parameters, specifically EP-SNR and PM-SNR at different noise levels and the standardization exponent, can significantly enhance localization performance, the study underscores the critical role of prior modeling in addressing depth bias. The incorporation of smoothing further refines the reconstructions, making the SKF approach competitive with more computationally intensive methods. Importantly, the use of SEP data provides a realistic and physiologically meaningful benchmark, reinforcing the relevance of the method for practical neuroimaging applications. These results contribute to the growing body of work that seeks to improve the spatial and temporal accuracy of non-invasive brain imaging and support the broader adoption of dynamic, Bayesian methods in clinical and research contexts.

This research contributes to the evolving field of EEG source localization by advancing the modeling of prior structures within a dynamic Bayesian framework for reconstructing deep brain activity. Although static methods such as MNE, sLORETA, and wMNE have been widely used, employing prior assumptions based on minimum norm, depth weighting, or anatomical constraints—they often lack the capacity to capture the temporal evolution of brain signals \citep{hamalainen1994interpreting, pascual2002standardized}. To address these limitations, dynamic approaches have been introduced, which model the temporal progression of neural sources through state-space formulations. Among these, Kalman filtering has emerged as a prominent method, offering recursive estimation and improved temporal resolution. It has been applied in various contexts for source localization over the years \citep{galka2004solution, long2006large}.

Different prior models have also been integrated into a dynamic framework for deep activity localization, for example, the SESAME approach uses sequential Monte Carlo techniques to infer dynamic dipole parameters non-parametrically, allowing for flexible modeling of sparse deep sources \citep{luria2024sesameeg}. RAMUS, on the other hand, incorporates anatomical and physiological constraints to enhance source recovery in deep regions by modeling recurrent cortical-subcortical dynamics \citep{rezaei2021reconstructing}. These methods demonstrate the diversity of prior formulations that can be used to support deep source estimation in dynamic settings, ranging from particle-based samplers to biologically informed priors. Within this broader context, our approach builds on the Gaussian prior formulation in the SKF, a recent extension of the Kalman framework that incorporates a depth-compensating standardization term into the recursive estimation process which offers a computationally efficient alternative that allows for explicit parameter tuning. In this work, we extend the idea of standardization into the smoothing stage, where basic RTS smoothing lacks inherent depth weighting and tends to attenuate source amplitudes. To compensate for this, we employed a higher standardization exponent of 1.25, building on findings from static methods, where MNE performs optimally with a depth exponent between 0.6 and 0.8 \citep{lin2006assessing}. This adjustment accounts for the bidirectional nature of smoothing and improves the reconstruction of deep sources.

This study contributes to the growing body of work on Gaussian prior models for deep brain activity reconstruction within dynamic Bayesian frameworks. Building on the SKF, which incorporates depth compensation into recursive estimation, we extend standardization into the smoothing phase to address the amplitude attenuation typical of basic RTS smoothing. By tuning the standardization exponent to 1.25, we demonstrate that Gaussian priors, when properly parameterized, can support accurate spatio-temporal deep source localization in EEG.

Despite its promising results, this study has several limitations. First, the evaluation was performed using simulated SEPs, while, offering controlled conditions and known ground truth, they do not fully capture the variability and complexity of real EEG data. The performance of SKF was assessed in real-world scenarios before. Thus, the results in this study remain to be validated experimentally. Second, while Gaussian priors offer computational tractability and interpretability, they may not fully represent the sparsity or non-Gaussian characteristics of actual brain activity, potentially limiting the reconstruction fidelity in some cases. Third, the computational burden of the SKF framework, particularly when extended with smoothing, restricts its scalability to very high-resolution source spaces or long-duration recordings without further optimization. Lastly, the process of selecting prior parameters still involves empirical tuning; although we adopted a data-driven approach, more principled or automated methods for prior selection would enhance reproducibility and applicability.

Future work could explore alternative prior structures beyond the Gaussian, such as sparse Bayesian learning or hierarchical models, to improve flexibility and capture a broader range of source dynamics. Adaptive standardization methods \citep{lahtinen2024standardized} and more sophisticated smoothing strategies could further enhance depth localization \citep{wan2000unscented, paul2009rssi}, especially in noisy conditions. Moreover, testing the approach on real EEG datasets, including clinical applications like epilepsy, would help validate the generalizability of the method. 

\section{Acknowledgements}

This study was supported by the Research Council of Finland (RCF) through the Center of Excellence in Inverse modeling and Imaging 2018--2025 (353089) and the  Flagship of Advanced Mathematics for Sensing, Imaging and modeling (FAME) (359185), Doctoral Education Pilot on Advanced Mathematics for Modelling, Sensing, and Imaging (DREAM), Ministry of Education and Culture, Finland, VN/3137/2024. The work of J.\ Lahtinen was supported by the Jenny and Antti Wihuri Foundation. For travel support, we are grateful to the joint DAAD/RCF researcher exchange project (RCF 367453), which has provided us the chance to visit Insitute for Biomagnetism and Biosignalanalysis (IBB), University of Münster, Münster, Germany. AK and SP were also supported by project 359198.

%%
%% End of file `elsarticle-template-num.tex'.

\bibliographystyle{elsarticle-num}
\bibliography{sn-bibliography}% common bib file

\begin{thebibliography}{10}
\expandafter\ifx\csname url\endcsname\relax
  \def\url#1{\texttt{#1}}\fi
\expandafter\ifx\csname urlprefix\endcsname\relax\def\urlprefix{URL }\fi
\expandafter\ifx\csname href\endcsname\relax
  \def\href#1#2{#2} \def\path#1{#1}\fi

\bibitem{laxton2013deep}
A.~W. Laxton, A.~M. Lozano, Deep brain stimulation for the treatment of alzheimer disease and dementias, World neurosurgery 80~(3-4) (2013) S28--e1.

\bibitem{obeso2008functional}
J.~A. Obeso, M.~C. Rodr{\'\i}guez-Oroz, B.~Benitez-Temino, F.~J. Blesa, J.~Guridi, C.~Marin, M.~Rodriguez, Functional organization of the basal ganglia: therapeutic implications for parkinson's disease, Movement disorders: official journal of the Movement Disorder Society 23~(S3) (2008) S548--S559.

\bibitem{seeber2019subcortical}
M.~Seeber, L.-M. Cantonas, M.~Hoevels, T.~Sesia, V.~Visser-Vandewalle, C.~M. Michel, Subcortical electrophysiological activity is detectable with high-density {EEG} source imaging, Nature communications 10~(1) (2019) 753.

\bibitem{galka2004solution}
A.~Galka, O.~Yamashita, T.~Ozaki, R.~Biscay, P.~Vald{\'e}s-Sosa, A solution to the dynamical inverse problem of {EEG} generation using spatiotemporal kalman filtering, NeuroImage 23~(2) (2004) 435--453.

\bibitem{somersalo2003non}
E.~Somersalo, A.~Voutilainen, J.~Kaipio, Non-stationary magnetoencephalography by {B}ayesian filtering of dipole models, Inverse Problems 19~(5) (2003) 1047.

\bibitem{schmitt2002efficient}
U.~Schmitt, A.~Louis, Efficient algorithms for the regularization of dynamic inverse problems: I. theory, Inverse Problems 18~(3) (2002) 645.

\bibitem{brooks1999inverse}
D.~H. Brooks, G.~F. Ahmad, R.~S. MacLeod, G.~M. Maratos, Inverse electrocardiography by simultaneous imposition of multiple constraints, IEEE Transactions on Biomedical Engineering 46~(1) (1999) 3--18.

\bibitem{dannhauer2013spatio}
M.~Dannhauer, E.~L{\"a}mmel, C.~H. Wolters, T.~R. Kn{\"o}sche, Spatio-temporal regularization in linear distributed source reconstruction from {EEG/MEG}: a critical evaluation, Brain topography 26 (2013) 229--246.

\bibitem{sorrentino2009dynamical}
A.~Sorrentino, L.~Parkkonen, A.~Pascarella, C.~Campi, M.~Piana, Dynamical meg source modeling with multi-target bayesian filtering, Human brain mapping 30~(6) (2009) 1911--1921.

\bibitem{sorrentino2010particle}
A.~Sorrentino, Particle filters for magnetoencephalography, Archives of Computational Methods in Engineering 17~(3) (2010) 213--251.

\bibitem{luria2024sesameeg}
G.~Luria, A.~Viani, A.~Pascarella, H.~Bornfleth, S.~Sommariva, A.~Sorrentino, The {SESAMEEG} package: a probabilistic tool for source localization and uncertainty quantification in {M/EEG}, Frontiers in Human Neuroscience 18 (2024) 1359753.

\bibitem{luria2020towards}
G.~Luria, D.~Duran, E.~Visani, D.~Rossi~Sebastiano, A.~Sorrentino, L.~Tassi, A.~Granvillano, S.~Franceschetti, F.~Panzica, Towards the automatic localization of the irritative zone through magnetic source imaging, Brain Topography 33 (2020) 651--663.

\bibitem{ribeiro2004kalman}
M.~I. Ribeiro, Kalman and extended {K}alman filters: Concept, derivation and properties, Institute for Systems and Robotics 43~(46) (2004) 3736--3741.

\bibitem{wan2000unscented}
E.~A. Wan, R.~Van Der~Merwe, The unscented {K}alman filter for nonlinear estimation, in: Proceedings of the IEEE 2000 adaptive systems for signal processing, communications, and control symposium (Cat. No. 00EX373), Ieee, 2000, pp. 153--158.

\bibitem{barton2008evaluating}
M.~J. Barton, P.~A. Robinson, S.~Kumar, A.~Galka, H.~F. Durrant-Whyte, J.~Guivant, T.~Ozaki, Evaluating the performance of {K}alman-filter-based {EEG} source localization, IEEE transactions on biomedical engineering 56~(1) (2008) 122--136.

\bibitem{hamid2021source}
L.~Hamid, N.~Habboush, P.~Stern, N.~Japaridze, {\"U}.~Aydin, C.~H. Wolters, J.~C. Claussen, U.~Heute, U.~Stephani, A.~Galka, et~al., Source imaging of deep-brain activity using the regional spatiotemporal {K}alman filter, Computer methods and programs in biomedicine 200 (2021) 105830.

\bibitem{Lahtinen2024SKF}
J.~Lahtinen, P.~Ronni, N.~P. Subramaniyam, A.~Koulouri, C.~Wolters, S.~Pursiainen, Standardized {K}alman filtering for dynamical source localization of concurrent subcortical and cortical brain activity, Clinical Neurophysiology 168 (2024) 15--24.

\bibitem{pascual2002standardized}
R.~D. Pascual-Marqui, et~al., Standardized low-resolution brain electromagnetic tomography (s{LORETA}): technical details, Methods Find Exp Clin Pharmacol 24~(Suppl D) (2002) 5--12.

\bibitem{rauch1965maximum}
H.~E. Rauch, F.~Tung, C.~T. Striebel, Maximum likelihood estimates of linear dynamic systems, AIAA journal 3~(8) (1965) 1445--1450.

\bibitem{sarkka2023bayesian}
S.~S{\"a}rkk{\"a}, L.~Svensson, Bayesian filtering and smoothing, Vol.~17, Cambridge university press, Cambridge, UK, 2023.

\bibitem{simon2006optimal}
D.~Simon, Optimal State Estimation: Kalman, H Infinity, and Nonlinear Approaches, John Wiley \& Sons, Hoboken, NJ, 2006.

\bibitem{lahtinen2024standardized}
J.~Lahtinen, A.~Koulouri, S.~Rampp, J.~Wellmer, C.~Wolters, S.~Pursiainen, Standardized hierarchical adaptive {L}p regression for noise robust focal epilepsy source reconstructions, Clinical Neurophysiology 159 (2024) 24--40.

\bibitem{grewal2014kalman}
M.~S. Grewal, A.~P. Andrews, Kalman filtering: Theory and Practice with {MATLAB}, John Wiley \& Sons, Hoboken, NJ, 2014.

\bibitem{rezaei2020parametrizing}
A.~Rezaei, M.~Antonakakis, M.~Piastra, C.~H. Wolters, S.~Pursiainen, Parametrizing the conditionally gaussian prior model for source localization with reference to the {P20/N20} component of median nerve {SEP/SEF}, Brain sciences 10~(12) (2020) 934.

\bibitem{mazziotta2001probabilistic}
J.~Mazziotta, A.~Toga, A.~Evans, P.~Fox, J.~Lancaster, K.~Zilles, R.~Woods, T.~Paus, G.~Simpson, B.~Pike, et~al., A probabilistic atlas and reference system for the human brain: International consortium for brain mapping ({ICBM}), Philosophical Transactions of the Royal Society of London. Series B: Biological Sciences 356~(1412) (2001) 1293--1322.

\bibitem{tadel2011brainstorm}
F.~Tadel, S.~Baillet, J.~C. Mosher, D.~Pantazis, R.~M. Leahy, Brainstorm: A user-friendly application for {MEG/EEG} analysis, Computational intelligence and neuroscience 2011~(1) (2011) 879716.

\bibitem{desikan2006automated}
R.~S. Desikan, F.~S{\'e}gonne, B.~Fischl, B.~T. Quinn, B.~C. Dickerson, D.~Blacker, R.~L. Buckner, A.~M. Dale, R.~P. Maguire, B.~T. Hyman, et~al., An automated labeling system for subdividing the human cerebral cortex on {MRI} scans into gyral based regions of interest, Neuroimage 31~(3) (2006) 968--980.

\bibitem{he2020zeffiro}
Q.~He, A.~Rezaei, S.~Pursiainen, Zeffiro user interface for electromagnetic brain imaging: a {GPU} accelerated {FEM} tool for forward and inverse computations in {Matlab}, Neuroinformatics 18 (2020) 237--250.

\bibitem{galaz2023multi}
F.~Galaz~Prieto, J.~Lahtinen, M.~Samavaki, S.~Pursiainen, Multi-compartment head modeling in {EEG}: Unstructured boundary-fitted tetra meshing with subcortical structures, {PLoS One} 18~(9) (2023) e0290715.

\bibitem{pursiainen2016electroencephalography}
S.~Pursiainen, J.~Vorwerk, C.~H. Wolters, Electroencephalography (eeg) forward modeling via {H(div)} finite element sources with focal interpolation, Physics in Medicine \& Biology 61~(24) (2016) 8502.

\bibitem{buchner1994source}
H.~Buchner, M.~Fuchs, H.~A. Wischmann, O.~D{\"o}ssel, I.~Ludwig, A.~Knepper, P.~Berg, Source analysis of median nerve and finger stimulated somatosensory evoked potentials: multichannel simultaneous recording of electric and magnetic fields combined with {3D-MR} tomography, Brain topography 6 (1994) 299--310.

\bibitem{rezaei2021reconstructing}
A.~Rezaei, J.~Lahtinen, F.~Neugebauer, M.~Antonakakis, M.~C. Piastra, A.~Koulouri, C.~H. Wolters, S.~Pursiainen, Reconstructing subcortical and cortical somatosensory activity via the {RAMUS} inverse source analysis technique using median nerve {SEP} data, Neuroimage 245 (2021) 118726.

\bibitem{hamalainen1994interpreting}
M.~S. H{\"a}m{\"a}l{\"a}inen, R.~J. Ilmoniemi, Interpreting magnetic fields of the brain: minimum norm estimates, Medical \& biological engineering \& computing 32 (1994) 35--42.

\bibitem{long2006large}
C.~J. Long, P.~L. Purdon, S.~Temereanca, N.~U. Desai, M.~Hamalainen, E.~N. Brown, Large scale kalman filtering solutions to the electrophysiological source localization problem-a {MEG} case study, in: 2006 International Conference of the IEEE Engineering in Medicine and Biology Society, IEEE, 2006, pp. 4532--4535.

\bibitem{lin2006assessing}
F.-H. Lin, T.~Witzel, S.~P. Ahlfors, S.~M. Stufflebeam, J.~W. Belliveau, M.~S. H{\"a}m{\"a}l{\"a}inen, Assessing and improving the spatial accuracy in {MEG} source localization by depth-weighted minimum-norm estimates, Neuroimage 31~(1) (2006) 160--171.

\bibitem{paul2009rssi}
A.~S. Paul, E.~A. Wan, Rssi-based indoor localization and tracking using sigma-point kalman smoothers, IEEE Journal of selected topics in signal processing 3~(5) (2009) 860--873.

\end{thebibliography}

\end{document}